\documentstyle[12pt]{article}
\def\appendix#1{
\addtocounter{section}{1}
\setcounter{equation}{0}
\renewcommand{\thesection}{\Alph{section}}
\section*{Appendix \thesection\protect\indent #1}
\addcontentsline{toc}{section}{Appendix \thesection\ \ \ #1}
}
\newcommand{\tr}[1]{\,{\rm tr}\,#1\,}

\def\be{\begin{equation}}
\def\la{\label}
\def\ee{\end{equation}}
\def\bea{\begin{eqnarray}}
\def\eea{\end{eqnarray}}
\def\eps{\varepsilon}
\newcommand{\iv}[3]{\bar{\phi}({#1})^{#2}_{#3}{}} 
\newcommand{\ov}[3]{\phi({#1})^{#2}_{#3}{}}     
\newcommand{\pw}[0]{\check{W}}                  
\begin{document}
\title{
\begin{flushright}
{\small q-alg/9612032}
\end{flushright}
\vspace{2cm}
$R$-matrix Quantization of the Elliptic 
Ruijsenaars--Schneider model. 
}
\author{G.E.~Arutyunov,
\thanks{Steklov Mathematical Institute,
Vavilov 42, GSP-1, 117966, Moscow, Russia; arut@class.mi.ras.ru}\\
L.O.~Chekhov
\thanks{Steklov Mathematical Institute,
Vavilov 42, GSP-1, 117966, Moscow, Russia; chekhov@class.mi.ras.ru}\\
and\\
S.A.~Frolov\thanks{Steklov Mathematical Institute,
Vavilov 42, GSP-1, 117966, Moscow, Russia; frolov@class.mi.ras.ru}  
}
\date {}
\maketitle
\begin{abstract}
It is shown that the classical $L$-operator algebra 
of the elliptic Ruijsenaars-Schneider model can be 
realized as a subalgebra of the algebra of functions on the cotangent 
bundle over the centrally extended current group in two 
dimensions. It is governed by two dynamical 
$r$ and $\bar{r}$-matrices satisfying a closed system of equations.
The corresponding quantum $R$ and $\overline{R}$-matrices are found
as solutions to quantum analogs of these equations.
We present the quantum $L$-operator algebra and show that the 
system of equations on $R$ and $\overline{R}$ arises
as the compatibility condition for this algebra. 
It turns out that the $R$-matrix is twist-equivalent to
the Felder elliptic $R^F$-matrix with $\overline{R}$ playing the role
of the twist. The simplest representation of the quantum 
$L$-operator algebra corresponding to the elliptic 
Ruijsenaars-Schneider model is obtained.
The connection of the quantum $L$-operator algebra to the 
fundamental relation $RLL=LLR$ with Belavin's elliptic $R$ matrix 
is established.
As a byproduct of our construction, we find a new 
$N$-parameter elliptic solution to the classical Yang-Baxter equation.
\end{abstract} 

\newpage
\section{Introduction}
\setcounter{equation}{0}
The appearance of classical dynamical $r$-matrices \cite{BV,AT} 
in the theory  of integrable many-body systems raises an interesting 
problem of their quantization. On this way one may hope to separate the 
variables explicitly. 

At present, the classical dynamical $r$-matrices are known for the
rational, trigonometric \cite{AT,ABT} and elliptic \cite{Skl,BS} 
Calogero--Moser (CM) systems, as well as for their relativistic 
generalizations -- the rational, trigonometric \cite{AR,Sur} and elliptic 
\cite{NKSR,Sur1} Ruijsenaars--Schneider (RS) systems \cite{R}. 
It is recognized that dynamical systems of the Calogero type  can
be naturally understood in the framework of the Hamiltonian reduction
procedure \cite{OP,GNe}. Moreover, the reduction
procedure provides an effective scheme to compute the 
corresponding dynamical $r$-matrices \cite{ABT,AM}. 

Depending explicitly on the phase variables,
the dynamical $r$-matrices do not satisfy a single
closed  equation of the Yang--Baxter type, that makes the 
problem of their quantization rather nontrivial.
In \cite{ABB}, the spin generalization of the 
Calogero--Sutherland system was 
quantized by using the particular solution \cite{Fel} of the 
Gervais--Neveu--Felder equation \cite{GN,Fel}
and in \cite{BBB}, it was interpreted in terms of quasi-Hopf algebras. 
This system is not integrable, but zero-weight representations of the 
quantum $L$-operator algebra admit a proper number of commuting 
integrals of motion. However, it seems to be important to find such a 
quantum $L$-operator algebra for the Calogero-type systems that 
possesses a sufficiently large abelian subalgebra.

Recently, an algebraic scheme for quantizing the rational 
RS model in the $R$-matrix formalism was proposed \cite{AFwe}.
We introduced a special parameterization of the cotangent bundle over
$GL(N,{\bf C})$. In new variables, the standard symplectic structure
was described by a classical (Frobenius) 
$r$-matrix and by a new dynamical $\bar{r}$-matrix. The classical
$L$-operator was introduced as a special matrix function on the
cotangent bundle. The Poisson algebra of $L$ inherited from 
the cotangent bundle coincided with the $L$-operator algebra
of the rational RS model. It is this reason why we called 
$L$ the $L$-operator. Quantizing 
the Poisson structure for $L$, we found the quantum 
$L$-operator algebra and constructed its particular representation 
corresponding to the rational RS system. This quantum algebra 
has a remarkable property, namely, it possesses a family of $N$ 
mutually commuting operators directly on the algebra level.

It is well-known that the elliptic RS model is the most general among 
the systems of the CM and RS types. 
In this paper, we are aimed to include
this model in our scheme. Recall that the classical $L$-operator 
algebra for the elliptic RS model can be obtained by means of the 
Poissonian \cite{AFM} or the Hamiltonian reduction schemes 
\cite{AFM1}. In the first scheme, the affine Heisenberg double is 
used as the initial phase space and in the second one, 
the cotangent bundle over the centrally extended current group in 
two dimensions is considered. Thus, the appropriate phase space we 
choose to deal with the elliptic RS model is the cotangent bundle 
$T^\star {\widehat{GL}}(N)(z,\bar z)$ over the centrally extended 
group ${\widehat{GL}}(N)(z,\bar z)$ of double loops. 
The application of our approach \cite{AFwe} is not straightforward
since one should work with the infinite-dimensional phase space and, 
therefore, the correct description of the Poisson structure on 
$T^\star {\widehat{GL}}(N)(z,\bar z)$ in the desired parameterization 
requires an intermediate regularization. 

Describe briefly the content of the paper and the results 
obtained. In the second section we start with describing 
the Poisson structure on $T^\star {\widehat{GL}}(N)(z,\bar z)$
that depends on two complex parameters $k$ and $\alpha$. Then
we parametrize $T^\star {\widehat{GL}}(N)(z,\bar z)$
in a special way. The Poisson structure in new variables
is ill defined due to the presence of singularities. 
To overcome this problem, we introduce an intermediate regularization.
Removing the regularization we find that only for $\alpha=1/N$
the resulting Poisson structure is well defined. The value
$\alpha=1/N$ corresponds, in fact, to the case where only the 
$SL(N)(z,\bar{z})$-subgroup is centrally extended. The corresponding
Poisson structure is described by two matrices 
$\bf r$ and $\bar{\bf r}$, which depend on $N$ dynamical variables $q_i$.
It follows from the Jacoby identity (see Appendix A) that $\bf r$ is 
an $N$-parameter elliptic solution to the classical Yang-Baxter 
equation (CYBE). It is worthwhile to mention that the main elliptic 
identities (see Appendix B) follow from the fulfillment of the 
CYBE for $\bf r$.  We expect that the 
matrix $\bf r$ is related to a special Frobenius subgroup in 
$GL(N)(z,\bar z)$ as it was in the finite-dimensional case 
\cite{AFwe}.  The Jacoby identity also implies a closed system of 
equations on $\bf r$ and $\bar{\bf r}$. 

We define a special matrix function $L$ on  
$T^\star {\widehat{GL}}(N)(z,\bar z)$. 
We call this function the ``$L$-operator" 
since as we show the Poisson algebra of $L$ inherited from 
$T^\star {\widehat{GL}}(N)(z,\bar z)$ literally 
coincides with the one for the elliptic RS model \cite{Sur1,AFM1}. 
Thus, the classical $L$-operator algebra can be 
realized as a subalgebra of the algebra of functions on 
$T^\star {\widehat{GL}}(N)(z,\bar z)$.
It turns out that  the $L$-operator 
as a function on $T^\star {\widehat{GL}}(N)(z,\bar z)$ admits
a factorization $L=WP$, where $p_i=\log P_i$ are the variables
canonically conjugated to $q_i$ and $W$ belongs to some special 
subgroup in $GL(N)(z,\bar z)$. The Poisson bracket for the entries 
of $W$ is given by the matrix $\bf r$ and coincides with the 
Sklyanin bracket defining the structure of the Poisson--Lie 
group.

Although the quantum analogs of equations on $\bf r$ and 
$\bar{\bf r}$ can be easily established, it is rather difficult to
find the corresponding quantum  $\bf R$ and 
$\overline{\bf R}$-matrices. The matter is that the matrices 
$\bf r$ and $\bar{\bf r}$ have the complicated structure,
${\bf r}=r-s\otimes I+I\otimes s $ and $\bar{\bf r}=\bar{r}-s\otimes I$,
due to the presence of the $s$-matrix. However, we observe that 
the classical $L$-operator algebra does not depend on $s$ and, 
moreover, the matrices $r$ and $\bar{r}$ also obey a closed
system of equations. We show that this system arises as the 
compatibility condition for a new Poisson algebra.
This algebra contains both the Poisson algebra of $T^\star 
{\widehat{GL}}(N)(z,\bar z)$  and the classical 
$L$-operator algebra as its subalgebras.  
In the third section, using this key observation, we pass 
to the quantization.

We find the corresponding quantum $R$ and $\overline{R}$-matrices
as solutions to quantum analogs of the equations for  $r$ and $\bar{r}$. 
In particular, the $R$-matrix satisfies a novel triangle relation
that differs from the standard quantum Yang-Baxter equation 
by shifting the spectral parameters in a special way.
The Felder elliptic $R^F$-matrix naturally arises in our 
construction. It turns out that the $R$-matrix is twist-equivalent to
the $R^F$-matrix with the $\overline{R}$-matrix playing the role
of the twist. 

Then we derive a new quadratic algebra satisfied by the ``quantum" 
$L$-operator. This algebra is described by the quantum dynamical
$R$-matrices, namely, $R$, $R^F$, and $\overline{R}$:
$$
R_{12}L_1\overline{R}_{21}L_2=L_2\overline{R}_{12}L_1R^{F}_{12}.
$$
We show that the system of equations on $R$, $R^F$, and 
$\overline{R}$-matrices arises as the compatibility condition for 
this algebra. We present the simplest representation of the quantum 
$L$-operator algebra corresponding to the elliptic RS model.
We note that when performing a simple canonical transformation, 
the quantum $L$-operator coincides in essential with the classical
$L$-operator found in \cite{R}.

The quantum integrals of motion for
the elliptic RS model were obtained in \cite{R}. In \cite{H},
it was shown that any operator from the Ruijsenaars commuting family 
can be realized as the trace of a proper transfer matrix for the
special $\hat{L}$-operator that obeys the relation 
$R\hat{L}\hat{L}=\hat{L}\hat{L}R$ with Belavin's elliptic $R$-matrix 
\cite{B}. We note that our $L$-operator is gauge-equivalent 
to $\hat{L}$. It follows from this observation that the determinant
formula for the commuting family \cite{H} is also valid for $L$.
We show that any representation of our $L$-operator
algebra is gauge equivalent to a representation of the relations
$R\hat{L}\hat{L}=\hat{L}\hat{L}R$. 

In Conclusion we discuss some problems to be solved.

\section {Classical $L$-operator algebra} 
\setcounter{equation}{0}
\subsection{Poisson structure of $T^\star {\widehat{GL}}(N)(z,\bar z)$} 

Let $T_\tau$ be a torus endowed with the standard complex structure and 
periods $1$ and $\tau$. Denote by $G$ a group of smooth mappings 
from $T_\tau$ into the group $GL(N,{\bf C})$. Then $g\in G$ is a 
double-periodic matrix function $g(z,\bar z)$. The dual space to the 
Lie algebra of $G$ is spanned by double-periodic functions 
$A(z,\bar z)$ with values in $\hbox{Mat\,}(N,{\bf C})$. In what 
follows, we often use the concise notation $g(z,\bar z)=g(z)$ and 
$A(z,\bar z)=A(z)$.  The group $G$ admits central extensions 
$\widehat{G}$~\cite{EF}.  The Poisson structure on $T^*\widehat{G}$ 
with fixed central charges reads
\begin{eqnarray} 
\label{AA}
\{ A_1(z),A_2(w)\} &=&\frac 1 
2 [C,A_1(z) -A_2(w)] \delta(z-w) -k (C-\alpha I)
\frac{\partial}{\partial \bar z }\delta(z-w) \\ 
\{g_1(z),g_2(w)\} &=& 0\label{gg}\\ 
\{A_1(z),g_2(w)\} &=&g_2 (w)C\delta (z-w), 
\label{pb} 
\end{eqnarray} 
where $k,\alpha$ are central charges and $\delta(z)$ is the 
two-dimensional $\delta$-function. Here we use the standard tensor 
notation, and $C$ is the permutation operator.  

One can consider the following Hamiltonian action of $G$ on 
$T^*\widehat{G}$ 
\begin{eqnarray} 
\label{gt} 
A(z)&\to 
&T^{-1}(z)A(z)T(z) +kT^{-1}(z)\bar\partial T(z),\\ \nonumber g(z)&\to 
&T^{-1}(z)g(z)T(z).  
\end{eqnarray} 
We restrict our consideration to 
the case of smooth elements $A(z)$. Then generic element $A(z)$ can 
be diagonalized by the transformation (\ref{gt})~\cite{FG}:  
\begin{equation} 
A(z)= T(z)DT^{-1}(z) 
-k\bar\partial T(z)T^{-1}(z).
\label{fac} 
\end{equation} 
Here $D$ is 
a constant diagonal matrix with entries $D_i$, $D_i\neq D_j$ and 
$T(z)$ is double-periodic.  Matrix $D$ is defined up to the action of 
the elliptic Weyl group.  One can fix $D$ by choosing the fundamental 
Weyl chamber.

Matrix $T(z)$ in eq.(\ref{fac}) is not uniquely defined.
Any element $\tilde{T}(z,{\bar z})=T(z,{\bar z})h(z)$, where
a diagonal matrix $h(z)$ is an entire function of $z$,
also satisfies (\ref{fac}). Demanding $\tilde{T}(z,{\bar z})$ to be 
double-periodic, we obtain that $h(z)$ is a constant matrix.
We can remove this ambiguity by imposing the condition
\begin{equation}
T(\varepsilon)e=e,
\label{Te}
\end{equation}
where $e$ is a vector such that $e_i=1$ $\forall i$, and $\varepsilon$ 
is an arbitrary point on $T_\tau$. In what follows, we denote the 
matrix $T(z)$ that solves eq.(\ref{fac}) and satisfies eq.(\ref{Te}) 
by $T^{\varepsilon}(z)$. Such matrices evidently form a group.

Now we try to rewrite the Poisson structure (\ref{AA}) in terms of 
variables $T$ and $D$. Since $D_i$ are $G$-invariant functions, they 
belong to the center of (\ref{AA}) and, therefore, it is enough to 
calculate the bracket $\{T^\eps(z),T^\eps(w)\}$. However,
the straightforward calculation reveals 
that this bracket is ill defined.
So, we begin with calculating
the bracket $\{T^{\varepsilon}(z),T^{\eta}(w)\}$, where 
$T^{\varepsilon}(z)$ and $T^{\eta}(w)$ satisfy (\ref{Te}) at 
different points $\varepsilon$ and $\eta$,
\begin{equation}
\{T_{ij}^{\varepsilon}(z),T^{\eta}_{kl}(w)\}=
\sum_{mnps}\int d^2z'd^2w'
\frac{\delta T_{ij}^{\varepsilon}(z)}{\delta A_{mn}(z')}
\frac{\delta T_{kl}^{\eta}(w)}{\delta A_{ps}(w')}
\{A_{mn}(z'),A_{ps}(w')\}.
\label{TT}
\end{equation}
To calculate the functional derivative 
$\frac{\delta T_{ij}^{\varepsilon}(z)}{\delta A_{mn}(z')}$, we 
consider the variation of (\ref{fac}):  
\begin{equation} 
X(z)=t(z)D-Dt(z)-k \bar\partial t(z)+d,
\label{sd}
\end{equation}
where $X(z)=T^{-1}(z)\delta A(z)T(z)$, $t(z)=T^{-1}(z)\delta T(z)$
and $d=\delta D$. 

First, from (\ref{sd}), we immediately obtain
\be
\frac{\delta D_{i}}{\delta A_{kl}(z)}
=\frac{1}{(\tau -\overline{\tau})}  
T_{ik}^{-1~\varepsilon}(z)T_{li}^{\varepsilon}(z).
\la{da}
\ee

Let us introduce the  function  $\Phi(z,s)$ of two complex variables
\begin{equation} 
\Phi(z,s)= 
\frac{\sigma(z+s)}
{\sigma(z)\sigma(s)}
e^{-2\zeta(\frac{1}{2})zs}e^{2\pi is\frac{z-\bar z}{\tau 
-\overline{\tau}}}.  \nonumber \end{equation} Here $\sigma (z)$ and 
$\zeta (z)$ are the Weierstrass $\sigma$- and $\zeta$-functions with 
periods equal to $1$ and $\tau$.  The  function  $\Phi(z,s)$ is 
the only double-periodic solution to the following equation:
\begin{equation} 
\bar\partial 
\Phi(z,s)+\frac{2\pi is}{\tau -{\overline\tau}} \Phi(z,s)=2\pi i\delta(z). 
\nonumber 
\end{equation} 
It is also convenient to define $\Phi(z,0)$ as follows:
$$ 
\Phi(z,0)=\lim_{\varepsilon \to 0}\left(\Phi(z,\varepsilon) - 
\frac{1}{\varepsilon}\right) =\zeta(z)-2\zeta(\frac{1}{2})z+2\pi i 
\frac{z-\bar{ z}}{\tau-\overline{\tau}}.
$$ 
This function solves the equation 
$$ 
\bar{\partial}\Phi(z,0)=2\pi i\delta(z) -\frac {2\pi 
i}{\tau-\overline{\tau}}.  
$$ 
We introduce the notation 
$q_{ij}\equiv q_i -q_j$, where $q_i=\frac{\tau -\overline\tau}{ 
2\pi ik}D_i$.

Using these functions, one can write the solution to (\ref{sd})
obeying the condition $t(\varepsilon)e=0$ \cite{AFM1}:
\begin{equation}
t(z)=\kappa\sum_{i,j}\int d^2 w ( \Phi(\varepsilon-w,q_{ij})
X_{ij}(w)E_{ii}-
\Phi(z-w,q_{ij})
X_{ij}(w)E_{ij} ).
\label{t}
\end{equation}
Hereafter, we denote $\frac1{2\pi ik}$ by $\kappa$.

Performing the variation of eq.(\ref{t}) with respect to $A_{mn}(w)$ one gets
\begin{eqnarray*}
\frac{\delta T_{ij}^{\varepsilon}(z)}{\delta A_{mn}(w)}
&=&\kappa  
\left( \sum_{k}\Phi(\varepsilon-w,q_{jk})
T_{ij}^{\varepsilon}(z)T_{jm}^{-1~\varepsilon}(w)T_{nk}^{\varepsilon}(w)
\right. \\
&-&\left.\sum_{k}\Phi(z-w,q_{kj})
T_{ik}^{\varepsilon}(z)T_{km}^{-1~\varepsilon}(w)T_{nj}^{\varepsilon}(w)
\right).
\end{eqnarray*}

To compute the bracket (\ref{TT}), one needs the following relation 
between $T^{\varepsilon}(z)$ and $T^{\eta}(z)$:
\begin{equation}
T^{\varepsilon}(z)=T^{\eta}(z)H^{\eta \varepsilon},
\label{TT'}
\end{equation}
where $H^{\eta \varepsilon}$ is a constant diagonal matrix.

By direct computation, one finds
\begin{equation}
\frac{1}{\kappa}\{T_1^{\varepsilon}(z),T_2^{\eta}(w)\}=
T_1^{\varepsilon}(z)T_2^{\eta}(w)(
H_1^{\varepsilon\eta}H_2^{\eta\varepsilon}r_{12}^{\varepsilon\eta}(z,w)-
\alpha f^{\varepsilon\eta}(z,w)), 
\label{TT''}
\end{equation}
where 
\begin{eqnarray}
r_{12}^{\varepsilon\eta}(z,w)&=&
\sum_{ij} \Phi(\varepsilon-\eta,q_{ij})E_{ii}\otimes E_{jj}
+\sum_{ij} \Phi(z-w,q_{ij})E_{ij}\otimes E_{ji}\\
&-&\sum_{ij} \Phi(z-\eta,q_{ij})E_{ij}\otimes E_{jj}
+\sum_{ij} \Phi(w-\varepsilon,q_{ij})E_{jj}\otimes E_{ij}
\label{rr}
\end{eqnarray}
and
$$
f^{\varepsilon\eta}(z,w)=
\Phi(\varepsilon-\eta,0)+\Phi(\omega-\varepsilon,0)+\Phi(z-\omega,0)
-\Phi(z-\eta,0).
$$      
The bracket (\ref{TT''}) has the $r$-matrix form with the $r$-matrix
depending not only on coordinates $q_i$ but also on the additional 
variables $H$.
          
In the limit $\eta\to\eps$, one encounters the singularity.
This shows that the variable~$T(z)$ is not a good candidate to describe
the Poisson structure~(\ref{AA}). However, one can use the freedom
to multiply $T(z)$ by any functional of $A$. So, we introduce a 
new variable 
\begin{equation}
{\bf T}^{\eps}(z)=T^{\eps}(z)(\det{ T^{\eps}(\eps)})^{\beta}.
\label{detT}
\end{equation}
We use $\det T^{\eps}(\eps)$ in the definition of ${\bf T}^{\eps}(z)$
in order to have the group structure for the new variables.

Using the Poisson bracket (\ref{TT''})  one immediately finds
\begin{eqnarray}
\nonumber
\frac{1}{\kappa}\{{\bf T}_1^{\varepsilon}(z),{\bf T}_2^{\eta}(w)\} &=&
{\bf T}_1^{\varepsilon}(z){\bf T}_2^{\eta}(w)\left(
H_1^{\varepsilon\eta}H_2^{\eta\varepsilon}r_{12}^{\varepsilon\eta}(z,w)-
\alpha f^{\varepsilon\eta}(z,w) \right. \\
\nonumber
&+&  
\beta I\otimes \tr{_3
H_3^{\varepsilon\eta}H_2^{\eta\varepsilon}
r_{32}^{\varepsilon\eta}(\eps,w)}
+ \beta \tr{_3
H_1^{\varepsilon\eta}H_3^{\eta\varepsilon}r_{13}^{\varepsilon\eta}(z,\eta)}
\otimes I \\
&+& \left.\beta^{2} \tr{_{34}
H_3^{\varepsilon\eta}H_4^{\eta\varepsilon}
r_{34}^{\varepsilon\eta}(\eps,\eta)}
I\otimes I \right) ,
\label{TTn}
\end{eqnarray}
since $f(\eps, w)=f(z, \eta)=0$. 

Now we are going to pass to the limit $\eta\to 
\eps$.\footnote{Without loss of generality we assume that $\eps$
and $\eta$ are real.} 
For this purpose one should take into account the following
behavior of $H^{\eps\eta}$ when $\eta$ goes to $\eps$:
$H^{\eps\eta}=1+ (\eps-\eta)h+o(\eps-\eta)$, where $h$ is a constant 
diagonal matrix being the functional of $A$. 
It turns out that there exists a unique choice for 
$\alpha$ and $\beta$, namely, $\alpha=1/N$,  $\beta=-1/N$, for which
the singularities cancel and there is no contribution from the matrix $h$.
In the limit $\eta\to \eps =0$, 
\footnote{The $\varepsilon$-dependence can be easily restored
by shifting both $z$ and $w$ by $\varepsilon$.}
for these values of $\alpha$ and 
$\beta$, one gets 
\be \la{TTF} \frac{1}{\kappa}\{{\bf T}_1(z),{\bf 
T}_2(w)\} = {\bf T}_1(z){\bf T}_2(w){\bf r}_{12}(z,w).  \ee 
Here the 
limiting ${\bf r}$-matrix is given by \be {\bf 
r}_{12}(z,w)=r_{12}(z,w)-s(z)\otimes I + I\otimes 
s(w)-\frac{1}{N}\Phi(z-w,0)I\otimes I, \la{rbf} \ee where \bea 
r(z,w)&=&
\sum_{i\neq j}\Phi(q_{ij})E_{ii}\otimes E_{jj}+
\sum_{ij}\Phi(z-w,q_{ij})E_{ij}\otimes E_{ji}\nonumber\\
&{}&-\sum_{ij}\Phi(z,q_{ij})E_{ij}\otimes E_{jj}+      
\sum_{ij}\Phi(w,q_{ij})E_{jj}\otimes E_{ij}
\la{r}
\eea
and 
\be
\la{s}
s(z)=\frac{1}{N}\sum_{ij}\left(
\Phi(q_{ij})E_{ii}-\Phi(z,q_{ij})E_{ij}\right).
\ee
Here we denote by the function $\Phi(q_{ij})$ the regular part
of $\Phi(\eps, q_{ij})$ for $\eps\to 0$:
\bea
\nonumber
\Phi(q_{ij})&=&\zeta(q_{ij})-2\zeta(\frac{1}{2})q_{ij}, ~~~i\neq 
j, \\ \nonumber \Phi(q_{ii})&=&0.  \eea Note that both ${\bf r}$ and 
$r$ are skew-symmetric:  ${\bf r}_{12}(z,w)=-{\bf r}_{21}(w,z)$.

A natural conjecture is that the ${\bf r}$-matrix obtained satisfies 
the classical Yang--Baxter equation
\be
[[{\bf r}, {\bf r}]]\equiv [{\bf r}_{12}(z_1,z_2), 
{\bf r}_{13}(z_1,z_3)+{\bf r}_{23}(z_2,z_3)]+
[{\bf r}_{13}(z_1,z_3), {\bf r}_{23}(z_2,z_3)]=0.
\la{cyb}
\ee
It can be verified either by direct calculation or by considering
the limiting case of the Jacoby identity for the bracket (\ref{TTn})
as is done in the Appendix A. Thereby, the ${\bf r}$-matrix (\ref{rbf})
is an $N$-parameter solution of the classical Yang--Baxter equation.

Let us note that, as one could expect, the condition $\det{\bf T}(0)=1$
is compatible with the bracket (\ref{TTF}), since $\det{\bf T}(z)$
is a central element of algebra (\ref{TTF}).

Remark that the choice $\alpha=1/N$ corresponds to the case where
only the $sl(N)(z,\bar{z})$-subalgebra is centrally extended. In terms of 
${\bf T}(z)$ ($\beta=-1/N$), the boundary condition looks like 
${\bf T}(0)e=\lambda e$ and $\det{{\bf T}(0)}=1$.
One can also check that the field $A(z)$ defined by (\ref{fac})
with the substitution ${\bf T}(z)$ for $T(z)$ obeys Poisson 
algebra (\ref{AA}).

The next step is to consider the special parameterization for the
field $g(z)$. To this end, we introduce $\tilde{A}(z)$:
\be
\tilde{A}=gAg^{-1}-k{\bar{\partial}}gg^{-1} +
\frac{k}{N}\tr{{\bar{\partial}}gg^{-1}}. 
\la{tilA}
\ee 
One can check that 
$\tilde{A}(z)$ Poisson commutes with $A(w)$ and obeys the Poisson 
algebra:  
\begin{eqnarray} \label{tilAA} \{ 
\tilde{A}_1(z),\tilde{A}_2(w)\} &=&-\frac 1 2 [C,\tilde{A}_1(z) 
-\tilde{A}_2(w)] \delta(z-w) +k (C-\frac{1}{N} I) 
\frac{\partial}{\partial \bar 
z }\delta(z-w) \nonumber\\ 
\{\tilde{A}_1(z),g_2(w)\} &=&Cg_2 (w)\delta (z-w). \label{tilpb} 
\end{eqnarray} 
Now we factorize $\tilde{A}(z)$  in the same manner as it was done
for $A(z)$,
\begin{equation} 
\tilde{A}(z)= {\bf U}(z)D{\bf U}^{-1}(z) 
-k\bar\partial {\bf U}(z){\bf U}^{-1}(z),
\label{fac'} 
\end{equation} 
where ${\bf U}(z)$ satisfies the boundary condition 
${\bf U}(0)e=\lambda e$ and $\det{{\bf U}(0)}=1$. Obviously,
${\bf U}(z)$ Poisson commutes with ${\bf T}(w)$ and satisfies the 
Poisson algebra
\be
\la{UUF}
\frac{1}{\kappa}\{{\bf U}_1(z),{\bf U}_2(w)\} =-{\bf U}_1(z){\bf U}_2(w){\bf r}_{12}(z,w).
\ee
One can find from (\ref{tilA}) and (\ref{fac'}) the representation
for the field $g$,
\be
g(z)=(\det g(z))^{\frac{1}{N}}{\bf U}(z){\bf P}{\bf T}^{-1}(z),
\la{g}
\ee
where ${\bf P}$ is a constant diagonal matrix. 

Computing the determinants  of the both sides of eq.(\ref{g}) one
gets 
$$\det{{\bf P}}=\det{({\bf T}(z)/{\bf U}(z))}.$$ 
Since the 
l.h.s. does not depend on $z$ and 
$\det{{\bf T}(0)}=\det{{\bf U}(0)}=1$ we obtain that $\det{{\bf P}}=1$ and 
$\det{{\bf T}(z)}=\det{{\bf U}(z)}$. 

Calculating the Poisson brackets of ${\bf P}$ with ${\bf P}$ and 
$Q=\mbox{diag}(q_1,\dots q_N)$ in the same manner as above one reveals that 
\bea 
\label{PQ}
\{{\bf P}_1,{\bf P}_2\}&=&0,\\
\frac{1}{\kappa}\{Q_1, {\bf P}_2\}&=&{\bf P}_2 (\sum_{ii}E_{ii}\otimes E_{ii}-\frac{1}{N}I\otimes I). 
\eea 
In fact, it means that $\log {\bf P}_i=p_i-\frac{1}{N}\sum_i p_i$, where
$p_i$ are canonically conjugated to $q_i$.

For the remaining Poisson brackets of ${\bf P}$ with the fields $\bf T$,  
$\bf U$, we have
\begin{eqnarray}
\label{TP}
\frac{1}{\kappa}\{{\bf T}_1(z), {\bf P}_2\}&=&{\bf T}_1(z){\bf P}_2\bar{{\bf r}}_{12}(z),\\
\label{UP}
\frac{1}{\kappa}\{{\bf U}_1(z), {\bf P}_2\}&=&{\bf U}_1(z){\bf P}_2\bar{{\bf r}}_{12}(z).
\end{eqnarray} 
Here 
\be
\bar{{\bf r}}_{12}(z)=\bar{r}_{12}(z)-
s(z)\otimes I-\frac{1}{N}I\otimes \sum_{ij} \Phi(q_{ij})E_{jj},
\la{bfchr}
\ee
where we introduced the $\bar{r}$-matrix:
\be
\bar{r}_{12}(z)=
\sum_{ij}\Phi(q_{ij})E_{ii}\otimes E_{jj}-
\sum_{ij}\Phi(z,q_{ij})E_{ij}\otimes E_{jj}.
\la{chr}
\ee

To complete the description of the classical Poisson structure of
the cotangent bundle we present the Poisson bracket of $\det g$ with
other variables:
\bea
\nonumber
\frac{1}{\kappa}\{Q,\det{g(w)}\}&=&\det{g(w)},\qquad
\{{\bf P},\det{g(w)}\}=0,\\
\nonumber
\frac{1}{\kappa}\{{\bf T}(z),\det{g(w)}\}&=&-\det{g(w)}{\bf 
T}(z)(\Phi(z-w,0)+\Phi(w,0)),\\
\nonumber
\frac{1}{\kappa}\{{\bf 
U}(z),\det{g(w)}\}&=&-\det{g(w)}{\bf U}(z)(\Phi(z-w,0)+\Phi(w,0)).
\eea

Recall that the Jacoby identity for bracket (\ref{TTF})
reduces to the classical Yang-Baxter equation for the 
${\bf r}$-matrix.  As to the Poisson relations (\ref{TP}) and (\ref{UP}), 
one finds that the Jacoby identity is equivalent to 
the following quadratic in $\bar{{\bf r}}$
equation:  
\be \la{jcr} [\bar{\bf 
r}_{12}(z),\bar{\bf r}_{13}(z)]-{\bf P}_3^{-1}\{\bar{\bf 
r}_{12}(z),{\bf P}_3\}+ {\bf P}_2^{-1}\{\bar{\bf r}_{13}(z),{\bf 
P}_2\}=0 
\ee 
and the equation involving ${\bf r}$ and $\bar{{\bf r}}$,
\be
\la{jrcr}
[{\bf r}_{12}(z,w),\bar{\bf r}_{13}(z)+\bar{\bf r}_{23}(w)]
+[\bar{\bf r}_{13}(z),\bar{\bf r}_{23}(w)]-
{\bf P}_3^{-1}\{\bar{\bf r}_{12}(z,w),{\bf P}_3\}=0. 
\ee
One can check by direct calculations that the matrices ${\bf r}$ and 
$\bar{\bf r}$ given by (\ref{rbf}) and (\ref{bfchr}) 
do solve these equations.

Finally, we remark that the fields $A(z)$ and $g(z)$ defined by 
(\ref{fac}) and (\ref{g}) obey  Poisson relations 
(\ref{AA}--\ref{pb}). 

Let us note that the Poisson relation for the generator 
${\bf W}(z)={\bf T}^{-1}(z){\bf U}(z)$ turns out to be the
Sklyanin bracket:
\be
\frac{1}{\kappa}\{{\bf W}_1(z),{\bf W}_2(w)\}=
[{\bf r}_{12}(z,w), {\bf W}_1(z){\bf W}_2(w)], 
\la{sWW} 
\ee
which, therefore, defines
the structure of a Poisson-Lie group. This group
is an infinite-dimensional analog of the Frobenius group appeared in
\cite{AFwe}, where the Poisson-Lie group structure was 
related to the existence of a non-degenerate two-cocycle on the
corresponding Lie algebra. It would be interesting to find a
similar interpretation in the infinite-dimensional case.

\subsection{Classical $L$-operator}
In this subsection, we define a special function on the cotangent
bundle, which we call the {\it classical $L$-operator}. 
The motivation to treat this function as the $L$-operator
is that the Poisson algebra of $L$ is equivalent to
the one found in \cite{Sur1} for the $L$-operator of the elliptic
RS model.

Denote by $L$ the following function
\be
L(z)={\bf T}^{-1}(z)g(z){\bf T}(z)=(\det g(z))^{\frac{1}{N}}{\bf T}^{-1}(z)
{\bf U}(z){\bf P}.
\la{L}
\ee
By using the formulas of the previous subsection, one can easily derive
the Poisson bracket of $L$ and $Q$: 
\be
\frac{1}{\kappa}\{Q_1,L_2(z)\}=L_2(z)\sum_i E_{ii}\otimes E_{ii},
\la{QL}
\ee
and the Poisson algebra of the $L$-operator,
\begin{eqnarray}
\label{LL}
\frac{1}{\kappa}\{L_1(z),L_2(w)\}&=&
{\bf r}_{12}(z,w)L_1(z)L_2(w) \\
\nonumber
&+&L_1(z)L_2(w)(\bar{\bf r}_{12}(z)-\bar{\bf r}_{21}(w)-{\bf r}_{12}(z,w))\\
\nonumber
&+&L_1(z)\bar{\bf r}_{21}(w)L_2(w)
-L_2(w)\bar{\bf r}_{12}(z)L_1(z).
\end{eqnarray}

Clearly, the generators $Q$ and $L$ form a Poisson subalgebra in the 
Poisson algebra of the cotangent bundle. An important feature of 
this subalgebra is that $I_n(z)=\tr{L^n(z)}$ form a set 
of mutually commuting variables.

Just as in the finite-dimensional case \cite{AFwe}, one can see from
(\ref{QL}) that the $L$-operator admits the following factorization:
$L(z)=W(z)P$, where $Q$ and ${\log P}$ are canonically conjugated
variables, the $W$-algebra coincides with (\ref{sWW}), and the 
bracket of $W$ and $P$ is 
\begin{equation} 
\frac{1}{\kappa}\{W_1(z),P_2\}=-P_2[\bar{\bf r}_{12}(z),W_1(z)].
\label{sWP}
\end{equation} 

In fact, everything what we need to quantize the $L$-operator algebra
(\ref{LL}) is prepared. The problem of quantization is reduced to 
finding the quantum ${\bf R}$ and $\overline{\bf R}$-matrices satisfying 
the quantum analogs of eqs.(\ref{cyb}), (\ref{jcr}), and (\ref{jrcr}),
\begin{eqnarray}
\label{RR}
{\bf R}_{12}(z_1,z_2){\bf R}_{21}(z_2,z_1)&=&1,\\
\label{RRR}
{\bf R}_{12}(z_1,z_2){\bf R}_{13}(z_1,z_3){\bf R}_{23}(z_2,z_3)&=
&{\bf R}_{23}(z_2,z_3){\bf R}_{13}(z_1,z_3){\bf R}_{12}(z_1,z_2),\\ \label{RRbRb}
{\bf R}_{12}(z_1,z_2)\overline{\bf R}_{13}(z_1)\overline{\bf R}_{23}(z_2)&=&
\overline{\bf R}_{23}(z_2)\overline{\bf R}_{13}(z_1)P_3^{-1}{\bf R}_{12}(z_1,z_2)P_3,\\ 
\label{RbRb}
\overline{\bf R}_{12}(z)P_2^{-1}\overline{\bf R}_{13}(z)P_2&=&
\overline{\bf R}_{13}(z)P_3^{-1}\overline{\bf R}_{12}(z)P_3.
\end{eqnarray}
These matrices are assumed to have the standard behavior near 
$\hbar =0$:  
\begin{eqnarray} \nonumber {\bf R}=1+\hbar{\bf r} + 
o(\hbar),\qquad \overline{\bf R}=1+\hbar \bar{\bf r} + o(\hbar),  
\end{eqnarray} 
where $\hbar$ is a quantization parameter. 

The problem formulated seems to be rather complicated due to the 
presence of the $s$-matrix in the classical ${\bf r}$ and 
$\bar{\bf r}$-matrices.
However, Poisson algebra (\ref{LL}) possesses
an important property allowing one to avoid the problem at hand.
Namely, the matrix $s(z)$ coming both in ${\bf r}$ and $\bar{\bf r}$ drops
out from the r.h.s. of (\ref{LL}). Thereby, eq.(\ref{LL}) can 
be eventually rewritten as 
\begin{eqnarray}
\label{LLf}
\frac{1}{\kappa}\{L_1(z),L_2(w)\}&=&
r_{12}(z,w)L_1(z)L_2(w) \\
\nonumber
&-&L_1(z)L_2(w)(r_{12}(z,w)+\bar{r}_{21}(w)- \bar{r}_{12}(z))\\
\nonumber
&+&L_1(z)\bar{r}_{21}(w)L_2(w)
-L_2(w)\bar{r}_{12}(z)L_1(z).
\end{eqnarray}
Moreover, if we denote by $r^F_{12}$ the sum
\be
r^F_{12}(z,w)=r_{12}(z,w)+\bar{r}_{21}(w)- \bar{r}_{12}(z),
\la{rF}
\ee
then using (\ref{r}) and (\ref{chr}) we obtain
\be
r^F_{12}(z-w)=-\sum_{ij}\Phi(q_{ij})E_{ii}\otimes E_{jj}+
\sum_{ij}\Phi(z-w,q_{ij})E_{ij}\otimes E_{ji}.
\la{rf}
\ee
In this expression one can recognize the 
elliptic solution to the classical Gervais--Neveu--Felder
equation \cite{GN,Fel}:
\bea
\la{cgnf}
&&[r_{12}^{F}(z_1-z_2), r_{13}^{F}(z_1-z_3)+ r_{23}^{F}(z_2-z_3)]+
[r_{13}^{F}(z_1-z_3), r_{23}^{F}(z_2-z_3)]\\
\nonumber
&&-P_3^{-1}\{r_{12}^{F}(z_1-z_2),P_3 \}
+P_2^{-1}\{r_{13}^{F}(z_1-z_3),P_2 \}
-P_1^{-1}\{r_{23}^{F}(z_2-z_3),P_1 \}=0.
\eea
In fact, $r^{F}$ emerges as the semiclassical limit of the quantum 
$R$-matrix found in \cite{Fel}. 

The absence of the $s$-matrix in the resulting $L$-operator 
algebra and the appearance of the $r^{F}$-matrix
show that there may exist a closed system
of equations involving only $r$- and $\bar{r}$-matrices in the 
classical case, and $R$- and $\overline{R}$-matrices in the quantum one.
In the next subsection we find the desired system of equations and 
describe a Poisson structure for which these equations ensure
the fulfillment of the Jacoby identity.

Note that the algebra (\ref{LLf}) literally coincides  with the one
obtained in \cite{AFM1} by using of the Hamiltonian reduction 
procedure.  A mere similarity transformation of $L$ turns 
algebra (\ref{LLf}) to the one previously found in \cite{Sur1}.  
In contrast to \cite{Sur1} where (\ref{LLf}) was derived
by direct calculation with the usage of the particular form
of the $L$-operator for the RS model, our treatment does not appeal 
to the particular form of $L$.

\subsection{Quadratic Poisson algebra with derivatives}
In the first subsection, we obtained the matrices ${\bf r}$ and
$\bar{\bf r}$ obeying system of equations (\ref{cyb}), 
(\ref{jcr}) and (\ref{jrcr}). Clearly, these equations are not 
satisfied when substituting $r$ and $\bar r$ for
${\bf r}$ and $\bar{\bf r}$.  
However, computing the l.h.s. of these equations after this
substitution we arrive at surprisingly simple result:  
\bea
\la{rder}
&&[r_{12}(z_1,z_2), r_{13}(z_1,z_3)+ r_{23}(z_2,z_3)]+
[r_{13}(z_1,z_3), r_{23}(z_2,z_3)]=\\
\nonumber
&&
-(\partial_1+\partial_2)r_{12}(z_1,z_2)
+(\partial_1+\partial_3)r_{13}(z_1,z_3)
-(\partial_2+\partial_3)r_{23}(z_2,z_3),
\eea
\be
[\bar{r}_{12}(z),\bar{r}_{13}(z)]-P_3^{-1}\{\bar{r}_{12}(z),P_3\}+
P_2^{-1}\{\bar{r}_{13}(z),P_2\}=-\partial 
(\bar{r}_{12}(z)-\bar{r}_{13}(z)),
\la{chrr} 
\ee
and
\bea
\nonumber
&&[r_{12}(z_1,z_2), \bar{r}_{13}(z_1)+ \bar{r}_{23}(z_2)]+
[\bar{r}_{13}(z_1), \bar{r}_{23}(z_2)]-P_3^{-1}\{r_{12}(z_1,z_2),P_3\} \\
\la{rchrd}
&&=
-(\partial_1+\partial_2)r_{12}(z_1,z_2)
+\partial_1\bar{r}_{13}(z_1)
-\partial_2\bar{r}_{23}(z_2).
\eea
Here $\partial=\frac{\partial}{\partial x}$, where $x=\mbox{Re}~z$.
Note that eqs.(\ref{jcr}) and (\ref{jrcr}) are formulated with the 
help of $\bf P$. However, since all the matrices depend only on the 
difference $q_{ij}=q_i-q_j$, we simply replace $\bf P$ by $P$.

Comparing eqs.(\ref{rder}--\ref{rchrd}) for $r$ and $\bar{r}$
with (\ref{cyb}), (\ref{jcr}), and (\ref{jrcr}) 
for ${\bf r}$ and $\bar{\bf r}$, we come to the 
conclusion that the $s(z)$-matrix coming in ${\bf r}$ and $\bar{\bf 
r}$ effectively plays the role of the derivative with respect to the 
spectral parameter.

It is worth mentioning that eqs.(\ref{rder})-(\ref{rchrd}) obtained
for $r$ and $\bar r$ can be rewritten in the same form as 
eqs.(\ref{cyb}), (\ref{jcr}), and (\ref{jrcr}) if we replace $r$ and 
$\bar r$ by $r_{12}-\partial_1+\partial_2$ and 
$\bar{r}_{12}-\partial_1$.  In particular, for (\ref{rder}), we have 
\be
\la{rer}
[r_{12}-\partial_1+\partial_2, r_{13}-\partial_1+\partial_3]+ 
[r_{13}-\partial_1+\partial_3, r_{23}-\partial_2+\partial_3]+ 
[r_{12}-\partial_1+\partial_2, r_{23}-\partial_2+\partial_3]=0.
\ee
Thus, $r_{12}-\partial_1+\partial_2$ is a matrix first-order
differential operator satisfying the standard classical 
Yang-Baxter equation. Using this fact we write down the
Poisson algebra generated by the fields $T(z)$, $U(z)$, $Q$ and $P$,
having eqs.(\ref{rder})-(\ref{rchrd}) as the consistency 
conditions:
\bea
\la{tt}
\frac{1}{\kappa}\{T_1(z),T_2(w)\} &=&
T_1(z)T_2(w)r_{12}(z,w)+T'_1T_2-T_1T'_2 ,\\
\la{tt1}
\frac{1}{\kappa}\{U_1(z),U_2(w)\} &=&-
U_1(z)U_2(w)r_{12}(z,w)-U'_1U_2+U_1U'_2 ,\\
\la{tt2}
\frac{1}{\kappa}\{T_1(z),P_2\}&=&P_2T_1(z)\bar{r}_{12}(z)+P_2T_1'(z),\\
\label{tt3}
\frac{1}{\kappa}\{U_1(z),P_2\}&=&P_2U_1(z)\bar{r}_{12}(z)+P_2U_1'(z),\\ 
\la{tt4}
\{Q_1,P_2\}&=&P_2\sum_i E_{ii}\otimes 
E_{ii},~~\{Q_1,T_2\}=\{Q_1,U_2\}=0,
\eea
where $T'=\partial T$. 

It is worth mentioning that the Poisson structure 
(\ref{tt}-\ref{tt4}) is not compatible with the boundary condition
$T(0)e=\lambda e$.

Let us note that there exists a Poisson subalgebra of Poisson 
algebra (\ref{tt}-\ref{tt4}), formed by the generators:
$$ A(z)=T(z)DT^{-1}(z)-k\bar{\partial} T(z)T^{-1}(z), 
~~~g(z)=U(z)PT^{-1}(z)$$ 
that coincides with the Poisson algebra of the 
cotangent bundle with the central charge $\alpha=0$.

Defining the $L$-operator as $L(z)=T^{-1}(z)g(z)T(z)=T^{-1}(z)U(z)P$,
we get for $L$ algebra (\ref{LLf}) obtained previously. As in
the previous subsection the commutativity of $I_n(z)$ follows again from
the one of $g(z)$. 

The main advantage of Poisson algebra 
(\ref{tt}-\ref{tt4}) is that it can be easily quantized.

\section{Quantization}
\setcounter{equation}{0}
\subsection{Quantum $R$-matrices}
In this section, following the ideology of the Quantum Inverse
Scattering Method \cite{F,KS},
we quantize the classical $r$ and $\bar{r}$-matrices
and derive the quantum $L$-operator algebra.

We start with quantization of the relations (\ref{rder})-(\ref{rchrd}).
Let $T(z)$, $U(z)$ be  matrix generating functions being the formal
Fourier series in variables $x$ and $y$:
$$
T(z)=\sum_{mn}T_{mn}e^{2\pi i(mx+ny)},
~~U(z)=\sum_{mn}U_{mn}e^{2\pi i(mx+ny)}, $$ where 
$z=x+\tau y$.

Denote by $\cal A$ a free associative unital algebra over the 
field $\bf C$ generated by matrix elements of the
Fourier modes of $T(z)$, $U(z)$, and by the entries of the
diagonal matrices $P$ and $Q$ modulo the relations 
\begin{eqnarray} 
\label{qTT} 
&&T_1(z)T_2(w-\hbar)=T_2(w)T_1(z-\hbar)R_{12}(-\hbar,z,w),\\
\label{qUU}
&&U_1(z)U_2(w+\hbar)=U_2(w)U_1(z+\hbar)R_{12}(\hbar,z,w),\\
\label{qTP}
&&T_1(z+\hbar)P_2\overline{R}_{12}(\hbar,z)=P_2T_1(z), \\
\label{qUP}
&&U_1(z+\hbar)P_2\overline{R}_{12}(\hbar,z)=P_2U_1(z), \\
\label{qQP}
&&[Q_1,P_2]=-\hbar P_2\sum_i E_{ii}\otimes E_{ii}, \\
\nonumber
&&
[T_1(z),U_2(w)]=[T_1(z),Q_2]=[U_1(z),Q_2]=[P_1,P_2]=[Q_1,Q_2]=0.
\end{eqnarray} 
Here $R(\hbar,z,w)$ and $\overline{R}_{12}(\hbar,z)$ are double-periodic
matrix functions of spectral parameters. These functions also depend 
on the coordinates $q_i$ and have the following
semiclassical behavior at $\hbar=0$:
\begin{equation}
R=1+\hbar r+ o(\hbar),~~\overline{R}=1+\hbar \bar{r}+ o(\hbar).
\label{ap}
\end{equation} 

The next step is to find the conditions on   
$R$ and $\overline{R}$ that ensure the consistency of 
the defining relations for $\cal A$.  In the sequel we often use
$R(z,w)$ as a shorthand notation for $R(\hbar,z,w)$.

First, we write down the compatibility condition for algebra 
(\ref{qTT}) or (\ref{qUU}), which  reduces to the Quantum Yang-Baxter 
equation with spectral parameters shifted by $\hbar$,  
\begin{equation} 
R_{12}(z,w)R_{13}(z-\hbar,s-\hbar)R_{23}(w,s)= 
R_{23}(w-\hbar,s-\hbar)R_{13}(z,s)R_{12}(z-\hbar,w-\hbar).
\label{QYB}
\end{equation} 

Analogously to the classical case, one can introduce the 
following matrix differential operator ${\cal R}(z,w)=
e^{\hbar\frac{\partial}{\partial w}}
R(z,w) 
e^{-\hbar\frac{\partial}{\partial z}}$
in terms of which eq.(\ref{QYB}) reads as the standard Quantum 
Yang-Baxter equation
\begin{equation}
{\cal R}_{12}(z,w){\cal R}_{13}(z,s){\cal R}_{23}(w,s)=
{\cal R}_{23}(w,s){\cal 
R}_{13}(z,s){\cal R}_{12}(z,w).  
\label{QYBD} 
\end{equation} 

Relation (\ref{qTT}) also requires 
the fulfillment of the ``unitarity" condition for $R$,
\begin{equation}
\label{unit}
R_{12}(z,w)R_{21}(w,z)=1.
\end{equation} 

Analogously, we find the following compatibility 
conditions for (\ref{qTP}):  
\begin{equation}
P_3^{-1}\overline{R}_{12}(z)P_3\overline{R}_{13}(z-\hbar)=
P_2^{-1}\overline{R}_{13}(z)P_2\overline{R}_{12}(z-\hbar)
\label{RRC}
\end{equation} 
and
\begin{equation}
P_3^{-1}R_{12}(z,w)P_3\overline{R}_{13}(z-\hbar)\overline{R}_{23}(w)=
\overline{R}_{23}(w-\hbar)\overline{R}_{13}(z)R_{12}(z-\hbar,w-\hbar).
\label{RRCRC}
\end{equation}  

Now taking into account (\ref{ap}) one can easily see that in the 
semiclassical limit
$$
-\frac{1}{\kappa}\{\cdot ,\cdot\}=\lim_{\hbar\to 0}  
\frac{1}{\hbar}[\cdot ,\cdot ] 
$$ 
relations (\ref{qTT}-\ref{qQP}) 
determine Poisson structure (\ref{tt}-\ref{tt4}), while 
eqs.(\ref{QYB}), (\ref{RRC}), and (\ref{RRCRC}) turn into 
(\ref{rder}), (\ref{chrr}), and (\ref{rchrd}), respectively, in order 
$\hbar^2$.  In the first order in $\hbar$, the unitarity condition 
(\ref{unit})  requires $r$ to be skew-symmetric.  
Hence, the algebra ${\cal A}$ with defining 
relations (\ref{qTT})-(\ref{qQP}), where $R$ and $\overline{R}$ are 
the solutions of (\ref{QYB}-\ref{RRCRC}) obeying (\ref{ap}), is a 
quantization of the Poisson structure (\ref{tt}-\ref{tt4}).

Now we are in a position to find the matrices $R$ and $\overline{R}$
explicitly. We start with the $R$-matrix for which we assume the 
following natural ansatz:
\begin{eqnarray}
\nonumber
&&fR(\hbar,z,w)= 
\sum_{ij}\Phi(\hbar_1,q_{ij}+\hbar_2)E_{ii}\otimes E_{jj}+
\sum_{ij}\Phi(z-w+\hbar_3,q_{ij}+\hbar_4)E_{ij}\otimes E_{ji}\\
\label{anR}
&&-
\sum_{ij}\Phi(z+\hbar_5,q_{ij}+\hbar_6)E_{ij}\otimes E_{jj}+      
\sum_{ij}\Phi(w+\hbar_7,q_{ij}+\hbar_8)E_{jj}\otimes E_{ij}.
\end{eqnarray}
This form is 
compatible with the structure of the classical $r$-matrix.
Here $\hbar_1,\dots,\hbar_8$ are arbitrary parameters that should
be specified by eqs.(\ref{QYB}) and (\ref{unit}), and 
$f$ is a scalar function that may depend only on $\hbar_i$ 
and spectral parameters. It turns out that the parameters $h_i$ 
are almost uniquely fixed by the unitarity condition (\ref{unit}).
Substituting (\ref{anR}) into (\ref{unit}) and using the
elliptic function identities we obtain
$$
\hbar_2=\hbar_3=\hbar_4=\hbar_6=\hbar_8=0,~~\hbar_5=\hbar_1+\hbar_7,
~~f^2(z,w)={\cal P}(\hbar_1)-{\cal P}(z-w),
$$
where ${\cal P}(z)$ is the Weierstrass ${\cal P}$-function.
\noindent 
Now it is a matter of 
direct calculation to check that eq.(\ref{QYB}) holds for 
$\hbar_1=\hbar$. The remaining parameter $\hbar_7$
is inessential since it corresponds to an
arbitrary common shift of the spectral parameters $z$ and $w$. 
In the sequel, we choose $\hbar_7=0$. Therefore, the obtained solution to
(\ref{QYB}) and (\ref{ap}) reads as follows:
\begin{eqnarray}
\label{QR}
f(z,w)R(\hbar,z,w)&=& 
\sum_{ij}\Phi(\hbar,q_{ij})E_{ii}\otimes E_{jj}+
\sum_{ij}\Phi(z-w,q_{ij})E_{ij}\otimes E_{ji}\\
&-&
\nonumber
\sum_{ij}\Phi(z+\hbar,q_{ij})E_{ij}\otimes E_{jj}+      
\sum_{ij}\Phi(w,q_{ij})E_{jj}\otimes E_{ij},
\end{eqnarray}
where $f(z,w)=\sqrt{{\cal P}(\hbar)-{\cal P}(z-w)}$.
One must be careful in the definition of $R(-\hbar,z,w)$. This matrix is
defined by (\ref{QR}) with the replacement $\hbar\to -\hbar$ and $f\to -f$.
Therefore, $R(\hbar,z,w)$ and $R(-\hbar,z,w)$ are related as 
\begin{equation}
R_{12}(-\hbar,z,w)=R_{21}(\hbar,w-\hbar,z-\hbar).
\label{min}
\end{equation} 

To find the $\overline{R}$-matrix, we adopt the following ansatz:
\begin{equation}
\label{anQRC}
\frac{1}{\sigma(\hbar)}\overline{R}_{12}(\hbar,z)=
\sum_{ij}\Phi(\hbar_1,q_{ij}+\hbar_2)E_{ii}\otimes E_{jj}-
\sum_{ij}\Phi(z+h_3,q_{ij}+\hbar_4 +\delta_{ij}\hbar_5 )E_{ij}\otimes 
E_{jj}.  \end{equation} It has almost the same matrix structure as the 
classical $\bar{r}$-matrix.  Since eq.(\ref{RRC}) is easier to deal 
with than eq.(\ref{RRCRC}), we first substitute (\ref{anQRC}) into 
eq.(\ref{RRC}) thus obtaining $\overline{R}$:  \begin{equation} 
\label{QRC}
\frac{1}{\sigma(\hbar)}\overline{R}_{12}(\hbar,z)=
\sum_{i\neq j}\Phi(\hbar,q_{ij})E_{ii}\otimes E_{jj}-
\sum_{i\neq j}\Phi(z+\hbar_{3} ,q_{ij})E_{ij}\otimes E_{jj}-
\Phi(z+\hbar_{3} ,-\hbar)\sum_{i}E_{ii}\otimes E_{ii}.
\end{equation} 
where $h_3$ remains unfixed.

Eq.(\ref{RRCRC}) involves both the $R$- and $\overline{R}$-matrices and
is independent on (\ref{QYB}) and (\ref{RRC}). One can 
verify by direct calculations that $R$ and $\overline{R}$
given by eqs.(\ref{QR}) and (\ref{QRC}) also satisfy (\ref{RRCRC})
as soon as $\hbar_3=\hbar$. 

One can easily check that in the case of real $\hbar$, the matrices
$R$ and $\overline{R}$ have the proper semiclassical behavior (\ref{ap}).

In what follows we also need the $\overline{R}^{-1}$-matrix,
\begin{equation}
\frac{1}{\sigma(\hbar)}\overline{R}_{12}^{-1}(\hbar,z)=
-\sum_{ij}\Phi(-\hbar,q_{ij}+\hbar)E_{ii}\otimes E_{jj}+
\sum_{ij}\Phi(z,q_{ij}+\hbar)E_{ij}\otimes E_{jj},
\label{inchR}
\end{equation} 

It would be of interest to mention that just as in the rational
case without the spectral parameter \cite{AFwe}, one can introduce 
the formal variable $W(z)=T^{-1}(z)U(z)$ with permutation relations
following from (\ref{qTT}--\ref{qQP}):
\begin{eqnarray}
\label{WW}
R_{12}(z,w)W_1(z)W_2(w+\hbar)&=&W_2(w)W_1(z+\hbar)R_{12}(z,w),\\
\label{WP}
W_1(z+\hbar )P_2\overline{R}_{12}(z)&=&  
P_2\overline{R}_{12}(z)W_1(z).  
\end{eqnarray} 
In analogy with the rational case, it is natural to treat eq.(\ref{WW}) 
as the defining relation of the quantum elliptic Frobenius group.

\subsection{Quantum $L$-operator algebra}
Just as in the classical case, we introduce a new variable:
\be
L(z)=T^{-1}(z)U(z)P=W(z)P,
\la{qL}
\ee
which we call a {\it quantum $L$-operator}.
Using the relations of the algebra $\cal A$ one can formally derive 
the following algebraic relations satisfied by the quantum $L$-operator:  
\begin{eqnarray}
\la{LQ}
[Q_1,L_2(z)]&=&-\hbar L_2(z)\sum_i E_{ii}\otimes E_{ii},\\
\nonumber
R_{12}(z,w)L_1(z)\overline{R}_{21}(w)L_2(w)&=&
L_2(w)\overline{R}_{12}(z)L_1(z)\overline{R}_{21}(w-\hbar )
R_{12}(z-\hbar ,w-\hbar )\overline{R}^{-1}_{12}(z-\hbar ).\\
\label{Lop}
\end{eqnarray} 
In spite of the fact that $L$ has the form $L(z)=W(z)P$, we can
not reconstruct from eqs.(\ref{LQ}) and (\ref{Lop}) the
relations (\ref{WW}) and (\ref{WP}) for $W$ and $P$. So, in 
the sequel, we do not assume any relations on $W$ and $P$.

Let us define 
\begin{equation}
R^F_{12}(z,w)=
\overline{R}_{21}(w)R_{12}(z,w)\overline{R}^{-1}_{12}(z).
\label{Rf}
\end{equation} 
Then, by using the explicit form of the $R$- and $\overline{R}$-matrices
and elliptic function identities, we obtain
\begin{equation}
fR^F_{12}(z-w)=
-\sum_{i\neq j}\Phi(-\hbar,q_{ij})E_{ii}\otimes E_{jj}
+\sum_{i\neq j}\Phi(z-w,q_{ij})E_{ij}\otimes E_{ji}+      
\Phi(z-w,\hbar)\sum_{i}E_{ii}\otimes E_{ii},              
\label{Rf1}
\end{equation} 
which is nothing but the $R$-matrix 
by Felder \cite{Fel}, i.e., an elliptic solution to 
the quantum Gervais--Neveu--Felder equation \cite{GN,Fel}:
\begin{equation}
P_1^{-1}R^{F}_{23}(w-s)P_1R^{F}_{13}(z-s)P_3^{-1}R^{F}_{12}(z-w)P_3=
R^{F}_{12}(z-w)P_2^{-1}R^{F}_{13}(z-s)P_2R^{F}_{23}(w-s).
\label{GNF}
\end{equation} 
Recall that one feature of $R^F$ is the ``weight zero" condition:
\begin{equation} 
[P_1P_2, R^{F}_{12}(z-w)]=0.
\label{wz}
\end{equation} 
Developing $R^F$ in powers of $\hbar$, we have
$R^F=1+\hbar r^F+ o(\hbar)$, where $r^F$ is given by (\ref{rf}).

Let us stress that in our consideration $R^F$ arises to account 
for the explicit form of $R$ and $\overline{R}$, and that the 
Gervais--Neveu--Felder equation does not follow 
from system (\ref{QYB}-\ref{RRCRC}).
Formula (\ref{Rf}) shows that the matrix $\overline{R}$ plays the role 
of the twist, which transforms the matrix $R(z,w)$ -- a particular 
solution of (\ref{QYB}) -- into such solution of (\ref{GNF}) which
depends only on the difference $z-w$. 

Thus, the quantum $L$-operator 
algebra (\ref{Lop}) can be presented in the following form:
\begin{equation}
R_{12}(z,w)L_1(z)\overline{R}_{21}(w)L_2(w)=
L_2(w)\overline{R}_{12}(z)L_1(z)R^F_{12}(z-w).
\label{LF}
\end{equation} 

The quantum $L$-operator algebra seems to be automatically 
compatible as $\cal A$ is compatible. However, 
a simple analysis shows that $\cal A$ and the algebra (\ref{LF})
admit different supplies of representations. In particular, the simplest
representation for $L$ we present below does not
realize the algebra (\ref{WW}), (\ref{WP}).
Therefore, we find it necessary to give a direct proof of the 
compatibility of (\ref{LF}). In this way, we come across 
eq.(\ref{GNF}) and discover a new relation involving $R^F$ and 
$\overline{R}$. To this end, let us multiply both sides of (\ref{LF}) 
by 
$$
P_2^{-1}\overline{R}_{31}(s+\hbar)P_2\overline{R}_{32}(s)L_3(s)$$ 
and subsequently using eq.(\ref{LF}) we transform the string 
$L_1\cdots L_2\cdots L_3$  
into $L_3\cdots L_2\cdots L_1$.  For the l.h.s., we have 
\begin{eqnarray} 
\nonumber && 
R_{12}(z,w)L_1\overline{R}_{21}(w)L_2
P_2^{-1}\overline{R}_{31}(s+\hbar)P_2\overline{R}_{32}(s)L_3= \\
\nonumber
&&
R_{12}(z,w)L_1\overline{R}_{21}(w)
\overline{R}_{31}(s+\hbar)L_2\overline{R}_{32}(s)L_3=\\
\nonumber            
&&
R_{12}(z,w)L_1\overline{R}_{21}(w)
\overline{R}_{31}(s+\hbar)R_{32}(s,w)L_3
\overline{R}_{23}(w)L_2R^{F}_{23}(w-s)=\\
\nonumber
&&       
R_{12}(z,w)R_{32}(s+\hbar,w+\hbar)
L_1\overline{R}_{31}(s)L_3P_3^{-1}\overline{R}_{21}(w+\hbar)P_3
\overline{R}_{23}(w)L_2R^{F}_{23}(w-s)=\\
\nonumber
&&
R_{12}(z,w)R_{32}(s+\hbar,w+\hbar)R_{31}(s,z)
L_3\overline{R}_{13}(z)L_1\times\\
\la{LLL}
&&~~~~~R^{F}_{13}(z-s)
P_3^{-1}\overline{R}_{21}(w+\hbar)P_3\overline{R}_{23}(w)L_2R^{F}_{23}(w-s).
\end{eqnarray}

At this point, we interrupt the chain of calculations by remarking 
that the next step implies the possibility to push $R^F$ somehow
through $P_3^{-1}\overline{R}_{21}(w+\hbar)P_3\overline{R}_{23}(w)$. 
It can be done by virtue of the following new relation
involving $R^F$ and $\overline{R}$:
\begin{equation}
R^F_{12}(z)P_2^{-1}\overline{R}_{31}(w)P_2\overline{R}_{32}(w-\hbar)=
P_1^{-1}\overline{R}_{32}(w)P_1\overline{R}_{31}(w-\hbar)R^F_{12}(z),
\label{RFCHR}
\end{equation} 
which can be checked directly by using the explicit forms (\ref{Rf1})  
and (\ref{QRC}) of $R^F$  and $\overline{R}$ respectively.

Now we pursue calculation (\ref{LLL}) with the relation 
(\ref{RFCHR}) at hand.  \begin{eqnarray} \nonumber && 
R_{12}(z,w)R_{32}(s+\hbar,w+\hbar)R_{31}(s,z)
L_3\overline{R}_{13}(z)\overline{R}_{23}(w+\hbar)\times \\
\nonumber
&&~~~~~~~~
L_1\overline{R}_{21}(w)L_2P_2^{-1}R^{F}_{13}(z-s)P_2R^{F}_{23}(w-s).
\end{eqnarray}
As to the r.h.s., the same technique yields
\begin{eqnarray}
\nonumber
&&
L_2\overline{R}_{12}(z)L_1R^{F}_{12}(z-w)           
P_2^{-1}\overline{R}_{31}(s+\hbar)P_2\overline{R}_{32}(s)L_3=\\
&&
\nonumber
L_2\overline{R}_{12}(z)\overline{R}_{32}(s+\hbar)L_1
\overline{R}_{31}(s)L_3P_3^{-1}R^{F}_{12}(z-w)P_3=\\
&&
\nonumber
L_2\overline{R}_{12}(z)\overline{R}_{32}(s+\hbar)
R_{31}(s,z)L_3\overline{R}_{13}(z)L_1
R^{F}_{13}(z-s)P_3^{-1}R^{F}_{12}(z-w)P_3\\
&&
\nonumber
=R_{31}(s+\hbar,z+\hbar)
L_2\overline{R}_{32}(s)L_3P_3^{-1}\overline{R}_{12}(z+\hbar)P_3 
\times \\
\nonumber
&&~~~~~~~~~
\overline{R}_{13}(z)L_1R^{F}_{13}(z-s)
P_3^{-1}R^{F}_{12}(z-w)P_3= \\ 
&& 
\nonumber 
R_{31}(s+\hbar,z+\hbar)R_{32}(s,w) 
L_3\overline{R}_{23}(w)\overline{R}_{13}(z+\hbar)L_2 
\overline{R}_{12}(z)L_1\times \\
&&~~~~~~~~~
P_1^{-1}R^{F}_{23}(w-s)P_1R^{F}_{13}(z-s)P_3^{-1}R^{F}_{12}(z-w)P_3=\\
&&
\nonumber
R_{31}(s+\hbar,z+\hbar)R_{32}(s,w)R_{12}(z+\hbar,w+\hbar)
L_3\overline{R}_{13}(z)\overline{R}_{23}(w+\hbar)L_1 
\times \\
\nonumber
&&~~~~~~~~~
\overline{R}_{21}(w)L_2R^{F}_{21}(w-z)
P_1^{-1}R^{F}_{23}(w-s)P_1R^{F}_{13}(z-s)P_3^{-1}R^{F}_{12}(z-w)P_3.
\end{eqnarray}
Therefore, comparing the resulting expressions we conclude that 
the compatibility condition for the $L$-operator algebra (\ref{LF})
reduces to four equations (\ref{QYB}), (\ref{RRCRC}), (\ref{Rf1}),
and (\ref{RFCHR}).

The existence of the Poisson commuting functions $I_n(z)$
in the classical case implies that the commuting family should exist 
in the quantum case as well. It should be the intrinsic property of 
the algebra (\ref{LF}) itself, without referring to the explicit form 
of its representations. Let us demonstrate the commutativity of the 
simplest quantities $\tr{L(z)}$ and $\tr{L^{-1}(z)}$ postponing the 
discussion of the general case to the next section.  To this end, we 
need one more relation involving the matrices $R^F$, $R$ and 
$\overline{R}$. 

In analogy with the rational case, it is useful
to introduce the variable $g(z)=U(z)PT^{-1}(z)$. 
Calculation of the commutator
$[g_1(z),g_2(w)]$ with the help of the defining relations of ${\cal A}$ 
results in 
\begin{eqnarray}
\nonumber
[g_1(z),g_2(w)]&=&U_2(w)U_1(z+\hbar)\left(
R_{12}(\hbar,z,w)P_1\overline{R}_{21}(w)P_2\overline{R}_{12}(z-\hbar) 
R_{12}(-\hbar,z,w)\right. \\
\nonumber
&-&\left. 
P_2\overline{R}_{12}(z)P_1 \overline{R}_{21}(w-\hbar)\right)
T_2^{-1}(w-\hbar)T_1^{-1}(z).
\end{eqnarray} 
When the spectral parameter is absent, the algebra $\cal A$ allows 
one to establish a connection with the quantum cotangent bundle (see \cite{AFwe}
for details). Then, in particular, the quantity  $[g_1,g_2]$ is equal 
to zero.  In the case at hand, we can not construct a subalgebra of
$\cal A$ that is isomorphic to the quantum cotangent bundle. However, 
one can note that in the elliptic case, the commutativity of $g(z)$ 
with $g(w)$ follows from the identity 
$$ 
R_{12}(\hbar,z,w)P_1\overline{R}_{21}(w)P_2\overline{R}_{12}(z-\hbar) 
R_{12}(-\hbar,z,w)=
P_2\overline{R}_{12}(z)P_1 \overline{R}_{21}(w-\hbar).
$$
Using the definition of $R^F$, eq.(\ref{min}),
and the "weight zero" condition (\ref{wz}), the last formula
can be written in the following elegant form
\begin{equation}
R_{12}(z,w)=P_2\overline{R}_{12}(z)P_2^{-1}
R^F_{12}(z-w)P_1\overline{R}^{-1}_{21}(w)P_1^{-1}.
\label{GG}
\end{equation} 
Identity (\ref{GG}) plays the primary role in proving the 
commutativity of the family $\tr{L(z)}$. To prove the commutativity, 
let us multiply both sides of 
$$
L_2(w)\overline{R}_{12}(z)L_1(z)=
R_{12}(z,w)L_1(z)\overline{R}_{21}(w)L_2(w)R^{F}_{21}(w-z).
$$
by $P_2\overline{R}^{-1}_{12}(z)P_2^{-1}$ and take the trace in the 
first and the second matrix spaces. We get \begin{eqnarray} \nonumber 
&&\tr{_{12}P_2\overline{R}^{-1}_{12}(z)P_2^{-1}
L_2(w)\overline{R}_{12}(z)L_1(z)}=\\
&&\tr{_{12}P_2\overline{R}^{-1}_{12}(z)P_2^{-1}
R_{12}(z,w)L_1(z)\overline{R}_{21}(w)L_2(w)R^{F}_{21}(w-z)}.
\label{lr}
\end{eqnarray}
 
It is useful to write $\overline{R}^{-1}_{12}$ in the factorized form
\begin{equation}
\overline{R}^{-1}_{12}(z)=
\sum_{ij}\overline{R}^{-1}_{ij}\otimes E_{jj}, 
\label{fchR} 
\end{equation} 
where 
$$
\overline{R}^{-1}_{ij}=-\sigma(\hbar)\Phi(-\hbar,q_{ij}+\hbar)E_{ii}+
\sigma(\hbar)\Phi(z,q_{ij}+\hbar)E_{ij}.
$$
Then the l.h.s. of (\ref{lr}) reads as
$$
\sum_{ij}\tr{_{12}(P_j\overline{R}^{-1}_{ij}P_j^{-1}\otimes E_{jj})
L_2\overline{R}_{12}(z)L_1}. 
$$
Taking into account that $\overline{R}_{12}$ is diagonal in the second
matrix space and using the cyclic property of the trace, we obtain
$$
\sum_{ij}\tr{_{12}(P_j\overline{R}^{-1}_{ij}P_j^{-1}\otimes I)
(I\otimes LE_{jj})\overline{R}_{12}(z)L_1}. 
$$
Since $L=WP$, where all entries of $W$ commutes with $q_i$,
we arrive at
\begin{eqnarray}
\nonumber
&&\sum_{ij}\tr{_{12}(P_j\overline{R}^{-1}_{ij}P_j^{-1}\otimes I)
(I\otimes WP_jE_{jj})\overline{R}_{12}(z)L_1}=\\
\nonumber
&&\sum_{ij}\tr{_{12}L_2(\overline{R}^{-1}_{ij}\otimes E_{jj})
\overline{R}_{12}(z)L_1}=\tr{L(w)}\tr{L(z)}.
\end{eqnarray}
As to the r.h.s. of (\ref{lr}), we use identity (\ref{GG})
to rewrite it as 
$$
\tr{_{12}R^{F}_{12}(z-w)P_1\overline{R}^{-1}_{21}(w)P_1^{-1}
L_1(z)\overline{R}_{21}(w)L_2(w)R^{F}_{21}(w-z)}.
$$
Having in mind that $\overline{R}_{21}$ is diagonal in the first 
matrix space and taking into account the property (\ref{wz}) one 
can easily see that under the trace sign, the matrix $R^F$ can be 
pushed to the right where it cancels with $R^{F}_{21}$. Therefore, 
we get 
\begin{equation} \tr{_{12}P_1\overline{R}^{-1}_{21}(w)P_1^{-1} 
L_1(z)\overline{R}_{21}(w)L_2(w)}.
\label{ftr}
\end{equation} 
Now applying to this expression the technique we used above for the l.h.s
of (\ref{lr}) we conclude that eq.(\ref{ftr}) is equal to
$\tr{L(z)}\tr{L(w)}$. Thus, we proved that $\tr{L(z)}$ commutes
with $\tr{L(w)}$. Quite analogously one can prove that
$\tr{L^{-1}(z)}$ commutes with $\tr{L^{-1}(w)}$ and with $\tr{L(w)}$. 

Now we give an example of the simplest representation of algebra
(\ref{LQ}) and (\ref{LF}) associated with the elliptic RS model.
Namely, the following $L$-operator satisfies algebra (\ref{LF}):
\begin{equation}
L(z)=\sum_{ij}\Phi(z,q_{ij}+\gamma)b_j P_j E_{ij},
\label{QuanL}
\end{equation} 
where 
\begin{equation}
b_j=        
\prod_{a\neq j}\Phi(\gamma,q_{aj}).
\label{b}
\end{equation} 
Here the parameter $\gamma$ is a coupling constant of the elliptic RS 
model. This can be checked by straightforward calculations. Some
comments are in order. The $L$-operator of the form (\ref{QuanL}) 
was already appeared in \cite{AFM1} as a result of the Hamiltonian reduction 
procedure applied to $T^\star {\widehat{GL}}(N)(z,\bar z)$. To guess
the explicit form of $b_j$ one should note that in the rational limit
algebra (\ref{LF}) tends to the one obtained in \cite{AFwe},
where the coefficients $b_j$ are found to be 
\be
b_j=    
\prod_{a\neq j}
\frac{q_{aj}+\gamma}{q_{aj}}.
\la{br'}
\ee
Therefore, it is natural to assume that the elliptic analog of (\ref{br'})
is given by (\ref{b}).

It is worthwhile to mention that $b_j$ are not uniquely 
defined since one can perform a canonical transformation of 
$(q,p)$-variables. In particular, the variables \begin{equation} 
\tilde{b}_j=        
\prod_{a\neq j}
\frac{\sigma(q_{aj}+\gamma)}{\sigma(\gamma)\sigma(q_{aj})}
\label{cb}
\end{equation} 
are related to $b_j$ by the canonical transformation
$q_i\to q_i$ and 
$$P_i\to  e^{\alpha\sum_{a}q_{ai}}P_i,$$
where 
$\alpha=2\zeta(1/2)\gamma - i\frac{\gamma-\bar{\gamma}}{\tau-{\bar{\tau}}}$. 

We call this $L$ the {\it quantum $L$-operator of the elliptic
RS model}. Indeed, taking the Hamiltonian to be $H=\tr{L(z)}$
one can see that the quantum canonical transformation of the form:  
\be 
P_i^{R}= 
\prod_{a\neq i} \left(\frac{\sigma(q_{ai}+\gamma)}{\sigma(q_{ai})}\right)^{1/2} 
P_i\prod_{a\neq i}
\left(\frac{\sigma(q_{ai})}{\sigma(q_{ai}-\gamma)}\right)^{1/2},
\la{qq}
\ee
where $P_i^{R}$ is the momentum in the Ruijsenaars Hamiltonian,
turn $H$ into the first integral $S_1$ from the Ruijsenaars 
commuting family \cite{R}. Moreover, after the canonical
transformation (\ref{qq}), the $L$-operator (\ref{QuanL})
coincides in essential with the classical $L$-operator of
the RS model.

The generating function for the commuting family
in terms of $L$ can be written as
\be
I(z,\mu)=:\det{(L(z)-\mu)}:=\sum_{k}(-\mu)^{N-k}I_k(z),
\la{cf}
\ee
where the normal ordering $::$ means that all 
momentum operators are pushed to the right.  

\noindent It follows from the results obtained in the next section.

\section{Connection to the fundamental relation $RLL=LLR$}
\setcounter{equation}{0}
In this section we establish a connection of the quantum 
$L$-operator algebra (\ref{LF}) with the fundamental relation
$RLL=LLR$. 

In \cite{H}, the operators from the Ruijsenaars commuting family 
were obtained by using a special representation $\hat{L}$ of the 
algebra 
\begin{equation} R^B_{12}(z-w) \hat{L}_1(z)\hat{L}_2(w) = 
\hat{L}_2(w)\hat{L}_1(z) R^B_{12}(z-w), 
\label{RLL} 
\end{equation} 
where $R^B(z)$ is Belavin's $R$-matrix being an elliptic solution
to the quantum Yang-Baxter equation \cite{B}. 
The explicit form of $R^B$ we use here can be found in \cite{RT}. 
For reader's convenience we recall a construction of $\hat{L}$ \cite{QF,H1}.

Denote by ${\bf h}^*$ the weight space for $sl_N({\bf C})$ that
can be realized in ${\bf C}^N$ with a basis $\epsilon_i$,
$<\epsilon_i,\epsilon_j>=\delta_{ij}$, as the orthogonal complement to
$\sum_{i=1}^{N} \epsilon_i$. Let $\bar{\epsilon}_k$ be 
the orthogonal projection of $\epsilon_k$:
$\bar{\epsilon}_k=\epsilon_k-\frac{1}{N}\sum_{i=1}^{N} \epsilon_i$.

For each $q \in \mbox{\bf h}^*$ one can introduce 
the intertwining vectors \cite{Bax,JMO}
\begin{equation}
\ov     {z}{q+\hbar \bar{\epsilon}_k}{q}_j
={\theta_j(\frac{z}{N}-<q,\bar{\epsilon}_k>)} /{i\eta(\tau)},
\label{inter}
\end{equation}
where
$$
\theta_{j}(z)=\sum_{n \in \frac{N}{2}-j+N{\bf Z}}^{}
{\exp 2\pi i\left[n (z+\frac{1}{2}) + \frac{n^2}{2N}\tau
\right]}
$$
and $\eta(\tau)= p^{1/24}\prod_{m=1}^\infty (1-p^m)$
is the Dedekind eta function with $p=\exp 2\pi i\tau$.

\noindent
Following \cite{QF} we denote by 
$\iv{z}{q+\hbar \bar{\epsilon}_k}{q}^j$
the entries of the matrix inverse to
$\ov{z}{q+\hbar \bar{\epsilon}_k}{q}_j$. Then the orthogonality
relations read as follows
\begin{equation}
\sum_{j=1}^n
{\bar{\phi}(z)}^{q +\hbar\bar{\epsilon}_k}_{q}{}^j
\ov{z}{q +\hbar\bar{\epsilon}_{k'}}{q}_j
=\delta_{kk'}
\qquad
\sum_{k=1}^n
\ov{z}{q+\hbar\bar{\epsilon_k}}{q} _j
{\bar{\phi}(z)}^{q+\hbar\bar{\epsilon_k}}_{q}{}^{j'}
=\delta_{jj'}.
\label{ort}
\end{equation}

In the sequel, the following formula \cite{H} will be of intensive use
\begin{eqnarray}
\la{iden}
\sum_{m=1}^n
\iv  {z}{q'+\hbar\bar{\epsilon}_j}{q'}^{m}
\ov     {z}{q+\hbar\bar{\epsilon}_{i}}{q}_{m}
=\frac{\theta(z+<q',\bar{\epsilon}_j>-<q,\bar{\epsilon}_{i}>)}{\theta(z)}
\prod_{j'\neq j}
\frac{\theta(<q',\bar{\epsilon}_{j'}>-<q,\bar{\epsilon}_{i}>)}
     {\theta(<q',\bar{\epsilon}_{j'}>-<q',\bar{\epsilon}_j>)}
\end{eqnarray}
Here $\theta(z)$ denotes the Jacoby $\theta$-function 
$$
\theta(z)=-\sum_n e^{2\pi i(z+\frac{1}{2})(n+\frac{1}{2})+i\pi 
\tau(n+\frac{1}{2})^2}=\theta'(0)e^{-\zeta(\frac{1}{2})z^2}\sigma(z).
$$

It is shown in \cite{H1,H2} that the $\hat{L}$-operator
\begin{equation}
\hat{L}_{ij}(z)=
\sum_{k=1}^{N}
\ov     {z+\gamma N}{q+\hbar\bar{\epsilon}_k}{q}_i
\iv     {z}{q+\hbar \bar{\epsilon}_k}{q}^j
e^{\hbar\bar{\epsilon}_k\frac{\partial}{\partial q_k}}, 
\label{LH} 
\end{equation}
acting on the space of functions on ${\bf h}^*$ satisfies 
relation (\ref{RLL}). This $\hat{L}$ is an $N\times N$ generalization  
of the $2\times 2$ Sklyanin $L$-operator \cite{S}.

The intertwining vectors $\ov{z}{q+\hbar \bar{\epsilon}_k}{q}_j$
coming in the definition of $\hat{L}$ relate the matrix $R^B$ with
the Boltzmann weights for the $A_{n-1}^{(1)}$ face model. 
Recall \cite{JMO} that the nonzero Boltzmann weights
depending on the spectral parameter $z$ are explicitly given by  
\bea
\la{weig}
&&\pw
\left[\begin{array}{ccc}
                        & q +\hbar\bar{\epsilon_{i}} & \\
                q &z& q + 2\hbar\bar{\epsilon_{i}} \\
                        & q +\hbar\bar{\epsilon_{i}} &
                        \\\end{array}\right
]
= \frac {\theta(z+{\hbar})}{\theta({\hbar})}, \\
\nonumber
&&\pw
\left[\begin{array}{ccc}
                & q +\hbar\bar{\epsilon_{i}} & \\
        q
        &z&q+\hbar(\bar{\epsilon_{i}}+\bar{\epsilon_{j}})\\
                & q +\hbar\bar{\epsilon_{i}} &
                \\\end{array}\right]
= \frac {\theta(-z+q_{ij})}{\theta(q_{ij})}
\; (i\neq j), \\
\nonumber
&&\pw
\left[\begin{array}{ccc}
                & q +\hbar\bar{\epsilon_{i}} & \\
        q &z&q + \hbar(\bar{\epsilon_{i}}
        +\bar{\epsilon_{j}})\\
                & q +\hbar\bar{\epsilon_{j}} &
                \\\end{array}\right]
= \frac {\theta(z)}{\theta({\hbar})}
        \frac {\theta({\hbar}+q_{ij})}{\theta(q_{ij})}
\; (i\neq j),
\eea
where $q_{ij}= <q, \bar{\epsilon_i} - \bar{\epsilon_j}>$.

The relation between $R^B$ and the face weights 
is given by 
\bea 
\la{Rint}
&&\sum_{i'j'}
R^B(z-w)^{i'j'}_{ij}
\ov     {z}{q+\hbar \bar{\epsilon}_k}{q}_{i'}
\ov     {w}{q+\hbar(\bar{\epsilon}_k+\bar{\epsilon}_m)}
        {q+\hbar \bar{\epsilon}_k}_{j'}
=       \\
\nonumber
&&\sum_{s}
\ov     {w}{q+\hbar \bar{\epsilon}_s}{q}_{j}
\ov     {z}{q+\hbar(\bar{\epsilon}_k+\bar{\epsilon}_m)}
        {q+\hbar \bar{\epsilon}_s}_{i}
~~\pw
\left[\begin{array}{ccc}
                        & q+\hbar \bar{\epsilon}_k & \\
                q       & z-w & q+\hbar(\bar{\epsilon}_k+\bar{\epsilon}_m)\\
                        & q+\hbar \bar{\epsilon}_s & \\
\end{array}\right]
\eea
\noindent

In what follows we use the concise notation
$$
W_{s}^{k}[k+m]=
\pw
\left[\begin{array}{ccc}
                        & q+\hbar \bar{\epsilon}_k & \\
                q       & z-w & q+\hbar(\bar{\epsilon}_k+\bar{\epsilon}_m)\\
                        & q+\hbar \bar{\epsilon}_s & \\
\end{array}\right]
$$
Then the dual relation to (\ref{Rint}) is
\begin{equation}
\sum_{i'j'}
\iv {w}{q+\hbar\bar{\epsilon}_k}{q}^{j'}
\iv     {z}{q+\hbar(\bar{\epsilon}_k+\bar{\epsilon}_m)}{q+\hbar\bar{\epsilon}_k}
        ^{i'}
R^B(z-w)^{ij}_{i'j'} =
\sum_{s} 
W_{k}^{s}[k+m]
\iv  {z}{q+\hbar\bar{\epsilon}_s}{q}^{i}
\iv  {w}{q+\hbar(\bar{\epsilon}_k+\bar{\epsilon}_m)}{q+\hbar\bar{\epsilon}_s}
        ^{j}
\label{dual}
\end{equation}

In \cite{H} another $L$-operator $\tilde{L}$ appeared. It is related 
to $\hat{L}$ in the following way:  
\bea 
\la{TL} 
\hat{L}_{ij}(z) 
\rightarrow \tilde{L}_{ij}(z) &=&\sum_{i'j'} 
\iv {z}{q+\hbar\bar{\epsilon}_{i}}{q}^{i'}
\ov     {z}{q+\hbar\bar{\epsilon}_{j}}{q}_{j'}
\hat{L}_{i'j'}(z)
\\
\nonumber
&=&
\frac{\theta(z+\gamma+q_{ij})}{\theta(z)}
\prod_{n\neq i}
\frac{\theta(\gamma+q_{nj})}{\theta(q_{ni})}~
e^{\hbar\bar{\epsilon}_j\frac{\partial}{\partial q_j}}.
\eea

In what follows we need to remove from the quantum $L$-operator algebra 
(\ref{LF}) the nonholomorphic dependence on the spectral and 
quantization parameters. This can be achived by considering the 
following transformation of the $L$-operator 
\be
L(z)\to e^{\alpha(z)Q}e^{-\beta Q}L(z)e^{\beta Q}e^{-\alpha(z)Q},
\ee
where $\alpha(z)$ is an arbitrary function of the spectral parameter
and $\beta$ is a complex number. Since the transformed $L$-operator 
also has the structure $WP$, then the following formula is 
valid 
\be
L_2(w)e^{\alpha(z)Q_1}=e^{\alpha(z)Q_1}L_2(w)e^{\hbar 
\alpha(z)r_0}, 
\la{form}
\ee
where the notation $r_0=\sum_iE_{ii}\otimes E_{ii}$ was used.

Recalling that the $L$-operator (\ref{QuanL})
satisfies the quantum $L$-operator algebra (\ref{LF}) 
and using eq.(\ref{form}) one can
easily establish the algebra satisfied by the transformed $L$:
$$
\check{R}_{12}(z,w)L_1(z)\check{\overline{R}}_{21}(w)L_2(w)=
L_2(w)\check{\overline{R}}_{12}(z)L_1(z)\check{R}^F_{12}(z-w),
$$ 
where the matrices $\check{R}$, $\check{\overline{R}}$ and $\check{R}^F$
are 
\begin{eqnarray}
\la{tr}
&&\check{R}_{12}(z,w)=
e^{-\alpha(z)Q_1-\alpha(w)Q_2+\beta Q_2}
R_{12}(z,w)
e^{\alpha(z)Q_1+\alpha(w)Q_2-\beta Q_1}, \\
\la{tr1}
&&\check{\overline{R}}_{12}(z)=
e^{\hbar \alpha(z)r_0
-\alpha(z)Q_1+\beta Q_2}
\overline{R}_{12}(z)e^{\alpha(z)Q_1-\beta Q_1},\\
\la{tr2}
&&\check{R}^F_{12}(z,w)=\\
\nonumber
&&
e^{-\alpha(z)Q_1-\alpha(w)Q_2
+\frac{\hbar}{2}(\alpha(w)-\alpha(z)) r_0+\beta Q_1}
R^F_{12}(z,w)e^{\alpha(z)Q_1+\alpha(w)Q_2
+\frac{\hbar}{2}(\alpha(w)-\alpha(z)) r_0-\beta Q_2}. 
\end{eqnarray}
Since the transformation in question keeps the
form of the quantum $L$-operator algebra intact, the transformed matrices
$\check{R}$, $\check{\overline{R}}$ and $\check{R}^F$ also
satisfy all the compatibility conditions. In particular, the transformation
(\ref{tr2}) defines another solution of (\ref{GNF}). For $\beta=0$ this
was observed in \cite{ABB}.

For the particular choice $\beta=2\pi i\frac{h-\bar{h}}{\tau-\bar{\tau}}$
and $\alpha(z)=2\pi i\frac{z-\bar{z}}{\tau-\bar{\tau}}+\beta$ we find 
that the matrices $\check{R}$, $\check{\overline{R}}$ and 
$\check{R}^F$ are given by the same formulae (\ref{QR}), (\ref{QRC}), 
(\ref{Rf1}) with the total change of  
$\Phi(z,s)=\frac{\theta(z+s)}{\theta(z)\theta(s)}
\theta'(0)e^{2\pi i \frac{z-\bar{z}}{\tau-\bar{\tau}}}$ 
by the meromorphic function
$\frac{\theta(z+s)}{\theta(z)\theta(s)}$. The transformation
with such a choice of $\alpha$ and $\beta$ transforms 
(up to an unessential multiplier) the
$L$-operator (\ref{QuanL}) into
\be L_{ij}(z)= 
\frac{\theta(z+q_{ij}+\gamma)}{\theta(z)\theta(q_{ij}+\gamma)}
\prod_{n\neq j}\frac{\theta(q_{nj}+\gamma)}{\theta(q_{nj})}{\bf P}_j 
\la{ourL}
\ee 
which is a quasi-periodic meromorphic matrix function of the spectral 
parameter:
$$
L(z+1)=L(z),~~~
L(z+\tau)=e^{-2\pi i(\gamma+\hbar)}e^{-2\pi i Q}L(z)e^{2\pi i Q}.
$$
We assume that the $L$-operator is of the form $W{\bf P}$, where
${\bf P}_i=e^{\hbar\bar{\epsilon}_i\frac{\partial}{\partial q_i}}$.
The replacement of $P$ by ${\bf P}$ preserves all the consistency
conditions because the $R$-matrices depend only on the difference 
$q_i-q_j$. Thus, eq.(\ref{LF})
with $R$-matrices defined via 
$\Phi(z,s)=\frac{\theta(z+s)}{\theta(z)\theta(s)}$ refers to 
the meromorphic version of the quantum $L$-operator algebra while
(\ref{ourL}) provides its particular meromorphic representation.
In what follows we use only this meromorphic version.

Comparing (\ref{ourL}) to (\ref{TL}) we can read off that $L$ 
and $\tilde L$ are related in the following way 
\be 
\la{rel} 
\tilde{L}_{ij}(z)=\frac{\prod_{n\neq j}\theta(q_{nj})}
{\prod_{n\neq i}\theta(q_{ni})}L_{ij}(z) 
\ee
Since the combined transformation (\ref{TL}), (\ref{rel}) from
$\hat{L}$ to $L$ depends only on $q$ we can conjecture that {\it any}
representation $L$ of the quantum $L$-operator algebra (\ref{LF}) is
gauge equivalent to some representation $\hat{L}$ of (\ref{RLL}) with
a gauge-equivalence defined as 
\be
\la{equiv}
\hat{L}_{ij}(z)=\sum_{i'j'}
\ov {z}{q+\hbar\bar{\epsilon}_{i'}}{q}_{i}
\iv {z}{q+\hbar\bar{\epsilon}_{j'}}{q}^{j}
\frac{\prod_{n\neq j'}\theta(q_{nj'})}
{\prod_{n\neq i'}\theta(q_{ni'})}L_{i'j'}(z). 
\ee

Now we are in a position to prove this conjecture.
Suppose $\hat{L}$ be an abstract $L$-operator satisfying algebra (\ref{RLL})
and introduce $\tilde{L}$ by eq.(\ref{TL}). Assume that $\tilde{L}$
has the structure $W{\bf P}$, where the entries of the diagonal
matrix ${\bf P}$ are
${\bf P}_i=e^{\hbar\bar{\epsilon}_i\frac{\partial}{\partial q_i}}$
and the entries of $W$ commute with $q_i$.
Then substituting $\hat{L}$ expressed via $\tilde{L}$ in (\ref{RLL}) and 
performing the straightforward calculation with the use of
(\ref{dual}) one finds an algebra satisfied by $\tilde{L}$:
\bea
\nonumber
%\la{TilL}
&&\sum_{\stackrel{i'j'}{abcd}}W_{k}^{a}[a+c]\delta_{a+c,i+k}
{\bar{\phi}(z)_{q}^{q+\hbar\bar{\epsilon}_a~i'}
\phi(z)_{q}^{q+\hbar\bar{\epsilon}_b} }_{i'}
\bar{\phi}(w)_{q+\hbar \bar{\epsilon}_a}
^{q+\hbar(\bar{\epsilon}_a+\bar{\epsilon}_c)~j'}
{\bar{\phi}(w)_{q+\hbar\bar{\epsilon}_j}
^{q+\hbar(\bar{\epsilon}_j+\bar{\epsilon}_d)} }_{j'}
\tilde{L}_{bj}(z)\tilde{L}_{dl}(w)=\\
\nonumber
&&\sum_{\stackrel{i'j'}{abcd}}W_{a}^{j}[j+l]\delta_{j+l,a+c}
\bar{\phi}(w)_{q}^{q+\hbar\bar{\epsilon}_k~i'}
{\phi(w)_{q}^{q+\hbar\bar{\epsilon}_d} }_{i'} 
\bar{\phi}(z)_{q+\hbar \bar{\epsilon}_k}
^{q+\hbar(\bar{\epsilon}_k+\bar{\epsilon}_i)~j'}
{\bar{\phi}(z)_{q+\hbar\bar{\epsilon}_c}
^{q+\hbar(\bar{\epsilon}_c+\bar{\epsilon}_b)} }_{j'}
\tilde{L}_{dc}(w)\tilde{L}_{ba}(z).
\eea
Performing the summation in $i'$ and $j'$ with the help of (\ref{iden})
we obtain
\bea
\nonumber
%\la{TsL}
&&\sum_{abc}W_{k}^{b}[b+a]\delta_{a+b,i+k}
\frac{\theta(w+q_{ac}+\hbar\delta_{ab}-\hbar\delta_{jc})}{\theta(w)}
\prod_{n\neq a}
\frac{\theta(q_{nc}+\hbar\delta_{nb}-\hbar\delta_{jc})}
{\theta(q_{na}+\hbar\delta_{nb}-\hbar\delta_{ab})}
\tilde{L}_{bj}(z)\tilde{L}_{cl}(w)=
\\
\nonumber
&&\sum_{abc}W_{a}^{j}[j+l]\delta_{j+l,a+c}
\frac{\theta(z+q_{ib}+\hbar\delta_{ik}-\hbar\delta_{bc})}{\theta(z)}
\prod_{n\neq i}
\frac{\theta(q_{nb}+\hbar\delta_{nk}-\hbar\delta_{bc})}
{\theta(q_{ni}+\hbar\delta_{nk}-\hbar\delta_{ik})}
\tilde{L}_{kc}(w)\tilde{L}_{ba}(z).
\eea
 
Let us introduce an operator $L$ by inverting (\ref{rel}). Substituting
this $L$ in the last formula, taking into account the nonzero components 
of the face weights and multiplying both sides by the function
$$
\frac{\prod_{n\neq k}\theta(q_{nk})}{\prod_{n\neq j}\theta(q_{nj})}
\frac{\prod_{n\neq i}\theta(q_{ni}+\hbar \delta_{nk}-\hbar \delta_{ik})}
{\prod_{n\neq l}\theta(q_{nl}+\hbar \delta_{nj}-\hbar \delta_{jl})}
$$
we finally arrive to the algebra  satisfied by $L$:
\bea
\la{LLH}
&&\sum_{s}W_{k}^{k}[2k]\delta_{i}^{k}
\frac{\theta(w+q_{ks}+\hbar-\hbar\delta_{js})}{\theta(w)}
\frac{\prod_{n\neq k}
\theta(q_{ns}-\hbar\delta_{js})}
{\prod_{n\neq s}\theta(q_{ns}+\hbar\delta_{nj}-\hbar\delta_{js})}
L_{kj}(z)L_{sl}(w)+\nonumber\\
&&\sum_{s}W_{k}^{k}[k+i]
\frac{\theta(w+q_{is}-\hbar\delta_{js})}{\theta(w)}
\frac{\prod_{n\neq i}
\theta(q_{ns}+\hbar\delta_{nk}-\hbar\delta_{js})}
{\prod_{n\neq s}\theta(q_{ns}+\hbar\delta_{nj}-\hbar\delta_{js})}
L_{kj}(z)L_{sl}(w)+\nonumber\\
&&\sum_{s}W_{k}^{i}[i+k]
\frac{\theta(q_{ki}+\hbar)}{\theta(q_{ki}-\hbar)}
\frac{\theta(w+q_{ks}-\hbar\delta_{js})}{\theta(w)}
\frac{\prod_{n\neq k}
\theta(q_{ns}+\hbar\delta_{ni}-\hbar\delta_{js})}
{\prod_{n\neq s}\theta(q_{ns}+\hbar\delta_{nj}-\hbar\delta_{js})}
L_{ij}(z)L_{sl}(w)=\nonumber\\
&&\sum_{s}W_{j}^{j}[j+l]
\left(\frac{\theta(q_{lj}+\hbar)}{\theta(q_{lj}-\hbar)}+2\delta_{lj}\right)
\frac{\theta(z+q_{is}+\hbar\delta_{ik}-\hbar \delta_{ls})}{\theta(z)}
\times\nonumber\\
&&\qquad\qquad\qquad\qquad\qquad\qquad\qquad\qquad\times
\frac{\prod_{n\neq i}
\theta(q_{ns}+\hbar\delta_{nk}-\hbar\delta_{ls})}
{\prod_{n\neq s}\theta(q_{ns}+\hbar\delta_{nl}-\hbar\delta_{ls})}
L_{kl}(w)L_{sj}(z)+\nonumber\\
&&
\sum_{s}W_{l}^{j}[j+l]
\frac{\theta(z+q_{is}+\hbar\delta_{ik}-\hbar\delta_{js})}{\theta(z)}
\frac{\prod_{n\neq i}
\theta(q_{ns}+\hbar\delta_{nk}-\hbar\delta_{js})}
{\prod_{n\neq s}\theta(q_{ns}+\hbar\delta_{nj}-\hbar\delta_{js})}
L_{kj}(w)L_{sl}(z).
\eea
Here in the second and the third lines $i\neq k$. The ratio
of products of theta-functions occurring in each term in (\ref{LLH})
allows one to take off the sum over $s$, e.g., when $i\neq k$,
we have
\bea
\nonumber
&&\frac{\prod_{n\neq i}
\theta(q_{ns}+\hbar\delta_{nk}-\hbar\delta_{js})}
{\prod_{n\neq s}\theta(q_{ns}+\hbar\delta_{nj}-\hbar\delta_{js})}=\\
\nonumber
&&
\delta_{is}\left(
\delta_{ij}\frac{\theta(q_{ki})}{\theta(q_{ki}-\hbar)}|_{j\neq k}
+\delta_{kj}(1-\delta_{ij})+
\frac{\theta(q_{ki}+\hbar)\theta(q_{ji})}{\theta(q_{ki})
\theta(q_{ji}+\hbar)}|_{{j\neq i}\atop{j\neq k}}
\right)+\\
\nonumber
&&
\delta_{ks}(1-\delta_{kj})\left(
\delta_{ij}\frac{\theta(\hbar)}{\theta(q_{ik}+\hbar)}
+
\frac{\theta(q_{jk})\theta(\hbar)}{\theta(q_{jk}+\hbar)
\theta(q_{ik})}|_{j\neq i}
\right)+\\
\nonumber
&&
\delta_{js}\frac{\theta(q_{kj})\theta(\hbar)}{\theta(q_{kj}-\hbar)
\theta(q_{ji}+\hbar)}|_{j\neq i}.
\eea

To compare (\ref{LLH}) to (\ref{LF}) we rewrite relation 
(\ref{LF}) with the help of eq.(\ref{GG}) in the following form:
\begin{equation}
R^F_{12}(z-w)P_1\overline{R}^{-1}_{21}(w)P_1^{-1}L_1(z)
\overline{R}_{21}(w)L_2(w)=
P_2\overline{R}^{-1}_{12}(z)P_2^{-1}
L_2(w)\overline{R}_{12}(z)L_1(z)R^{F}_{12}(z-w).
\label{LFF}
\end{equation} 
In the component form algebra (\ref{LFF}) is
presented in Appendix C. Comparing the components of (\ref{LLH})
to the ones of (\ref{LFF}) we establish that they coincide
up to the overall multiplicative factor $\theta(z-w)\theta(\hbar)^2$. 
Thus, we have shown that any representation of algebra (\ref{LF})
by the transformation (\ref{equiv}) turns into a representation
of (\ref{RLL}).  The connection established gives right to assert 
that algebra (\ref{LF}) possesses a family of $N$-commuting integrals
and that the formula of the determinant type (\ref{cf}) for the 
commuting family proved in \cite{H} is also valid for the $L$-operator
(\ref{QuanL}).

\section{Conclusion}
\setcounter{equation}{0}
In this paper, we described the dynamical $R$-matrix structure 
of the quantum elliptic RS model. The quantum $L$-operator algebra 
possesses a family of commuting operators. It turns out that
this algebra has a surprisingly simple structure and can be analysed 
explicitly in the component form. Furthermore, one can hope 
that the problem of finding new representations of the algebra 
obtained is simpler than the corresponding problem for the algebra 
$RLL=LLR$.  

There are several interesting problems to be discussed. 

First, we recall that in the classical case we obtained 
two different Poisson algebras, which lead to the same classical 
$L$-operator algebra. Only one of them was quantized. It is desirable to 
quantize the second one and to show that the corresponding quantum 
$L$-operator algebra is isomorphic to the algebra 
obtained in the paper.

The elliptic RS model we dealt with in the paper corresponds to
the $A_{N-1}$ root system. It seems to be possible to extend our
approach to other root systems and to derive the corresponding 
$L$-operator algebras. To this end, one should find a proper
parameterization of the corresponding cotangent bundle. 

Generalizing  our approach to the cotangent bundle over
a centrally extended group of smooth mappings from a higher-genus Riemann 
surface into a Lie group, one may expect to obtain new integrable 
systems. 

It is known that the CM systems admit spin generalizations 
\cite{KBBT,Ne,Rub}. 
Recently, the spin generalization was found for the elliptic RS model 
\cite{KrZ}. However, the Hamiltonian formulation for the model
is not found yet. One may hope that in our approach the spin models 
can arise as higher representations of the $L$-operator algebra.

Probably, the most interesting and complicated problem is to
separate variables for the quantum elliptic RS model. 
Up to now only the three-particle case for the trigonometric RS model
was solved explicitly \cite{SklKu}.
One could expect that the $L$-operator algebra obtained in the
paper may shed light on the problem.

{\bf ACKNOWLEDGMENT} The authors are grateful to 
M.A.Olshanetsky and N.A.Slavnov for the valuable 
discussions.  This work is supported in part by the RFBI grants 
N96-01-00608, N96-02-19085 and N96-01-00551 and by the ISF grant 
a96-1516.  

\setcounter{section}{0}

\appendix{}
\setcounter{equation}{0}
In this Appendix, we prove that the limiting ${\bf r}$-matrix
(\ref{rbf}) satisfies the classical Yang--Baxter equation.
For this, let us write the equation, which follows from the Jacoby 
identity for the bracket (\ref{TTn}):
\bea
\nonumber
&&[{\bf r}^{\eps\eta}_{12}(z,w), {\bf r}^{\eps\rho}_{13}(z,s)+
{\bf r}^{\eta\rho}_{23}(w,s)]+ 
[{\bf r}^{\eps\rho}_{13}(z,s), {\bf r}^{\eta\rho}_{23}(w,s)]\\
\la{Jac}
&&
+\frac{1}{\kappa}T_3^{\rho~-1}(s)\{T_3^{\rho}(s),{\bf 
r}^{\eps\eta}_{12}(z,w)]+ 
\frac{1}{\kappa}T_1^{\eps~-1}(z)
\{T_1^{\eps}(z),{\bf r}^{\eta\rho}_{23}(w,s)]\\ 
\nonumber
&&- 
\frac{1}{\kappa}T_2^{\eta~-1}(w)\{T_2^{\eta}(w),{\bf r}^{\eps\rho}_{13}(z,s)]=0.
\eea
Here 
\begin{eqnarray}
\nonumber
{\bf r}^{\eps\eta}_{12}(z,w) &=&
H_1^{\varepsilon\eta}H_2^{\eta\varepsilon}r_{12}^{\varepsilon\eta}(z,w)-
\alpha f^{\varepsilon\eta}(z,w) \\
\nonumber
&+&  
\beta I\otimes \tr{_3
H_3^{\varepsilon\eta}H_2^{\eta\varepsilon}r_{32}^{\varepsilon\eta}(z,w)}
+ \beta \tr{_3
H_1^{\varepsilon\eta}H_3^{\eta\varepsilon}r_{13}^{\varepsilon\eta}(z,w)}
\otimes I \\
&+& \beta^{2} \tr{_{34}
H_3^{\varepsilon\eta}H_4^{\eta\varepsilon}r_{34}^{\varepsilon\eta}(z,w)}
I\otimes I.
\label{rad}
\end{eqnarray}
To study (\ref{Jac}), one needs to know the Poisson bracket of 
${\bf T}^{\eps}(z)$ with $H^{\eta\rho}$. This bracket can be easily 
derived from eqs.(\ref{TT'}), (\ref{TT''}), and (\ref{detT}):
\bea
\la{TH}
&&\frac{1}{\kappa}{\bf T}^{\eps~-1}_1(z)
\{{\bf T}^\eps_1(z), H^{\eta\rho}_2\} 
= 
\sum_{i\neq j}H^{\eps\rho}_iH^{\eta\eps}_j
\Phi(\eps-\rho,q_{ij}) (E_{ii}-\frac{1}{N}I)\otimes E_{jj} \\
\nonumber
&&-
\sum_{i\neq j}H^{\eps\eta}_iH^{\eta\eps}_jH^{\eta\rho}_j
\Phi(\eps-\eta,q_{ij})(E_{ii}-\frac{1}{N}I)\otimes E_{jj}   \\
\nonumber
&&-
\sum_{i\neq j}H^{\eps\rho}_iH^{\eta\eps}_j
\Phi(z-\rho,q_{ij}) E_{ij}\otimes E_{jj}+
\sum_{i\neq j}H^{\eps\eta}_iH^{\eta\eps}_jH^{\eta\rho}_j
\Phi(z-\eta,q_{ij})E_{ij}\otimes E_{jj}   \\
\nonumber
&&+
\sum_j\frac{1}{N}H^{\eta\rho}_j\left(
\Phi(\eps-\rho,0)+\Phi(\eta-\eps,0)+\Phi(z-\eta,0)
-\Phi(z-\rho,0)\right)(E_{jj}-\frac{1}{N}I)\otimes E_{jj}.
\eea
Eq.(\ref{Jac}) holds at arbitrary values of all the parameters.
Without loss of generality one can put $\rho=0$ and $\eta=a\eps$,
where $a$ and $\eps$ are real. Let us now
perform the change of variables $H_{i}^{\eps\eta}$:
$$
H_{i}^{\eps\eta}=1+h_{i}^{\eps\eta}
$$
and consider the expansion of the l.h.s. of (\ref{Jac}) in powers of
$\eps $ and $h$.

Let us note that the matrix ${\bf r}^{\eps\eta}(z,w)$  has the 
following expansion in powers of $h$ at $\eps-\eta\to 0$:  
\be
{\bf r}^{\eps\eta}(z,w)={\bf r}(z,w)+
{\bf r}_{\mbox{reg}}^{(1)}(z,w,\eps)h+
{\bf r}^{(2)}(z,w,\eps)h^2 +o(\eps-\eta),
\la{rex}
\ee
where ${\bf r}(z,w)$ is given by (\ref{rbf}). The matrix 
${\bf r}_{\mbox{reg}}^{(1)}(z,w,\eps)$ is regular when $\eps\to 0$
and the $\eps$-dependence of ${\bf r}^{(2)}(z,w,\eps)$ is inessential.
Substituting (\ref{rex}) into the 
bracket ${\bf T}^{-1}\{{\bf T},{\bf r}\}$ one gets:
\be
\la{Tr}
{\bf T}^{-1}\{{\bf T},{\bf r}\} =
{\bf r}_{\mbox{reg}}^{(1)}(\eps){\bf T}^{-1}\{{\bf T},h\}+
{\bf r}^{(2)}(\eps)h{\bf T}^{-1}\{{\bf T},h\}.
\ee
Clearly, eq.(\ref{Jac}) should be satisfied in any order in $h$ and $\eps$.
Since we are interested in finding an equation for ${\bf r}$, we consider 
only the terms independent of $h$ and $\eps$ in the expansion of 
eq.(\ref{Jac}). The terms of zero order in $h$ and $\eps$ occurring in 
(\ref{Jac}) come from ${\bf r}^{\eps\eta}$ and from the first term in 
(\ref{Tr}). However, from the explicit expression for the bracket 
$\{{\bf T},H\}$ one can see that
$$
\{{\bf T},h\}|_{h=0}=o(\eps)
$$
Now, taking into account that  ${\bf r}_{\mbox{reg}}^{(1)}(\eps)$ is 
regular at $\eps\to 0$, we conclude that the last three terms in 
(\ref{Jac}) do not contribute. Thus, in the zero order in $h$ and $\eps$, 
eq.(\ref{Jac}) reduces to the CYBE for ${\bf r}(z,w)$. 

The remarkable thing is that the main elliptic identities 
(see Appendix  B) follow from the Jacoby identity for the bracket (\ref{TTn})             
or, equivalently, from the Yang-Baxter equation for ${\bf r}(z,w)$.  
  
\appendix{} 
Here we present the basic elliptic function identities, 
formulated as a set of functional relations on $\Phi(z,w)$ \cite{NKSR}:  
\begin{eqnarray} 
\label{third} 
&&\Phi(z,x)\Phi(w,y)=\Phi(z,x-y)\Phi(z+w,y)+\Phi(z+w,x)\Phi(w,y-x),\\
\label{lim}
&&\Phi(z,x)\Phi(z,y)=\Phi(z,x+y)\left(
\Phi(z,0)+\Phi(x,0)+\Phi(y,0)-\Phi(z+x+y,0)\right), \\
\nonumber
&&\Phi(z-w,a-b)\Phi(z,x+b)\Phi(w,y+a)-
\Phi(z-w,x-y)\Phi(z,y+a)\Phi(w,x+b)=\\
&&\Phi(z,x+a)\Phi(w,y+b)
\left(\Phi(a-b,0)+\Phi(x+b,0)-\Phi(x-y,0)-\Phi(a+y,0)\right).
\label{cub}
\end{eqnarray} 
Eq.(\ref{lim}) is the limiting case of eq.(\ref{third}) where $w\to z$, 
and eq.(\ref{cub}) is a consequence of (\ref{third}) and (\ref{lim}). 
Note that the exponent term
$e^{-2\zeta(1/2)zs+2\pi is\frac{z-\bar{z}}{\tau-\overline{\tau}}}$
in $\Phi(z,s)$ as well as the linear term 
$-2\zeta(1/2)z+2\pi i\frac{z-\bar{z}}{\tau-\overline{\tau}}$ in
$\Phi(z,0)$ are irrelevant since they drop out from 
(\ref{third}--\ref{lim}).

To establish the unitarity relation for $R$, one
also needs the identity involving the Weierstrass ${\cal P}$-function
$$
\Phi(z,s)\Phi(z,-s)={\cal P}(z)-{\cal P}(s),
$$
and to prove eq.(\ref{rder}), the following relation between the 
derivatives of $\Phi$ is of use:
$$
\frac{\partial \Phi(z,q_{ij})}{\partial z}=
\frac{\partial \Phi(z,q_{ij})}{\partial q_{ij}}-
(\Phi(z,0)-\Phi(q_{ij}))\Phi(z,q_{ij}). 
$$
\appendix{}
In this Appendix we present the quantum $L$-operator algebra
\begin{equation}
R^F_{12}(z-w)P_1\overline{R}^{-1}_{21}(w)P_1^{-1}L_1(z)
\overline{R}_{21}(w)L_2(w)=
P_2\overline{R}^{-1}_{12}(z)P_2^{-1}
L_2(w)\overline{R}_{12}(z)L_1(z)R^{F}_{12}(z-w)
\label{LFF'}
\end{equation} 
in the component form. 
The l.h.s. of (\ref{LFF'}) has the form
\begin{eqnarray*}
&&\sum_{i\neq k;j,l}
\Phi(\hbar,q_{ki})\Phi(\hbar,q_{ik})\Phi(\hbar,q_{kj}-\hbar)L_{ij}(z)L_{kl}(w) 
E_{ij}\otimes E_{kl}\\
-&&\sum_{i\neq k;j,l}
\Phi(\hbar,q_{ki})\Phi(\hbar,q_{ij}-\hbar)\Phi(w,q_{kj}-\hbar)L_{ij}(z)L_{jl}(w) 
E_{ij}\otimes E_{kl}\\
+&&\sum_{i\neq k;j,l}
\Phi(\hbar,q_{ki})\Phi(w,q_{ki})\Phi(\hbar,q_{ij}-\hbar)L_{ij}(z)L_{il}(w)
E_{ij}\otimes E_{kl}\\
+&&\sum_{i\neq k;j,l}
\Phi(z-w,q_{ik})\Phi(\hbar,q_{ki})\Phi(\hbar,q_{ij}-\hbar)L_{kj}(z)L_{il}(w)
E_{ij}\otimes E_{kl}\\
-&&\sum_{i\neq k;j,l}
\Phi(z-w,q_{ik})\Phi(\hbar,q_{kj}-\hbar)\Phi(w,q_{ij}-\hbar)L_{kj}(z)L_{jl}(w)
E_{ij}\otimes E_{kl}\\
+&&\sum_{i\neq k;j,l}
\Phi(z-w,q_{ik})\Phi(w,q_{ik})\Phi(\hbar,q_{kj}-\hbar)L_{kj}(z)L_{kl}(w)
E_{ij}\otimes E_{kl}\\
+&&\sum_{j,k,l}
\Phi(z-w,\hbar)\Phi(w,\hbar)\Phi(\hbar,q_{kj}-\hbar)L_{kj}(z)L_{kl}(w)
E_{kj}\otimes E_{kl}\\
-&&\sum_{j,k,l}
\Phi(z-w,\hbar)\Phi(w,\hbar)\Phi(w+\hbar,q_{kj}-\hbar)L_{kj}(z)L_{jl}(w)
E_{kj}\otimes E_{kl}
\end{eqnarray*}
The r.h.s. of (\ref{LFF'}) reads 
\begin{eqnarray*}
&&\sum_{i\neq k;j\neq l}
\Phi(\hbar,q_{ki})\Phi(\hbar,q_{il}-\hbar)\Phi(\hbar,q_{lj})
L_{kl}(w)L_{ij}(z) E_{ij}\otimes E_{kl}\\
+&&\sum_{i\neq k;j,l}
\Phi(\hbar,q_{ki})\Phi(\hbar,q_{ij}-\hbar)
\Phi(z-w,q_{lj}+\hbar \delta_{lj})
L_{kj}(w)L_{il}(z) E_{ij}\otimes E_{kl}\\
-&&\sum_{i\neq k;j\neq l}
\Phi(\hbar,q_{kl}-\hbar)\Phi(z,q_{il}-\hbar)
\Phi(\hbar,q_{lj})
L_{kl}(w)L_{lj}(z) E_{ij}\otimes E_{kl}\\
-&&\sum_{i\neq k;j,l}
\Phi(\hbar,q_{kj}-\hbar)\Phi(z,q_{ij}-\hbar)
\Phi(z-w,q_{lj}+\hbar \delta_{lj})
L_{kj}(w)L_{jl}(z) E_{ij}\otimes E_{kl}\\
+&&\sum_{i\neq k;j\neq l}
\Phi(z,q_{ik})\Phi(\hbar,q_{kl}-\hbar)\Phi(\hbar,q_{lj})
L_{kl}(w)L_{kj}(z) E_{ij}\otimes E_{kl}\\
+&&\sum_{i\neq k;j,l}
\Phi(z,q_{ik})\Phi(\hbar,q_{kj}-\hbar)\Phi(z-w,q_{lj}+\hbar 
\delta_{lj})L_{kj}(w)L_{kl}(z) E_{ij}\otimes E_{kl}\\
+&&\sum_{j\neq l;i}
\Phi(z,\hbar)\Phi(\hbar,q_{il}-\hbar)\Phi(\hbar,q_{lj})
L_{il}(w)L_{ij}(z) E_{ij}\otimes E_{il}\\
+&&\sum_{i,j,l}
\Phi(z,\hbar)\Phi(\hbar,q_{ij}-\hbar)\Phi(z-w,q_{lj}+\hbar 
\delta_{lj})L_{ij}(w)L_{il}(z) E_{ij}\otimes E_{il}\\
-&&\sum_{j\neq l;i}
\Phi(z,\hbar)\Phi(z+\hbar,q_{il}-\hbar)\Phi(\hbar,q_{lj})
L_{il}(w)L_{lj}(z) E_{ij}\otimes E_{il}\\
-&&\sum_{i,j,l}
\Phi(z,\hbar)\Phi(z+\hbar,q_{ij}-\hbar)\Phi(z-w,q_{lj}+\hbar 
\delta_{jl}) L_{ij}(w)L_{jl}(z) E_{ij}\otimes E_{il}
\end{eqnarray*}
%%%%%%%%%%%%%%%%%%%%%%%%%%%%%%%%%%%%%%%%%%%%%%%%%%%%%%%%%%%%%%%%

\end{document}